\documentclass[preprint,12pt]{elsarticle}
\usepackage{amsmath}
\usepackage{amssymb}
\usepackage{lscape}
\usepackage{capt-of}
\usepackage{float,lineno}
\usepackage{color,colortbl}
\usepackage{natbib}
\usepackage{graphicx}
\newcommand{\R}{\mathbb{R}}
\newcommand{\myd}{\text{ d}}
\journal{ }
\begin{document}
\begin{frontmatter}
%
%
%\title{Comparative sensitivity analysis of muscle activation dynamics}
\title{Comparing different muscle activation dynamics using sensitivity analysis }
\author[koblenz]{Robert Rockenfeller \corref{corauthor}}
\ead{rrockenfeller@uni-koblenz.de}
\author[stgt1,jena]{Michael G\"unther}
\author[stgt1,stgt2]{Syn Schmitt}
\author[koblenz]{Thomas G\"otz}
\address[koblenz]{
Universit\"at Koblenz,
Institut f\"ur Mathematik,
56070 Koblenz,
Deutschland}
\address[stgt1]{
Universit\"at Stuttgart,
Institut f\"ur Sport- und Bewegungswissenschaft,
Allmandring 28,
70569 Stuttgart,
Deutschland}
\address[jena]{
Friedrich--Schiller--Universit\"at,
Institut f\"ur Sportwissenschaft,
Lehrstuhl f\"ur Bewegungswissenschaft,
Seidelstra{\ss}e 20,
07749 Jena,
Deutschland}
\address[stgt2]{
Stuttgart Research Centre for Simulation Technology,
Pfaffenwaldring 7a,
D-70569 Stuttgart,
Deutschland}
\begin{abstract} 
In this paper, we mathematically compared two models of mammalian
striated muscle activation dynamics proposed by \citet{Hatze} and
\citet{Zajac}. Both models are representative of a broad variety of
biomechanical models formulated as ordinary differential equations
(ODEs). The models incorporate some parameters that directly represent
known physiological properties. Other parameters have been introduced
to reproduce empirical observations. We used sensitivity analysis as a
mathematical tool for investigating the influence of model parameters
on the solution of the ODEs. That is, we adopted a former approach
\citep{Lehman1982a} for calculating such (first order)
sensitivities. Additionally, we expanded it to treating initial
conditions as parameters and to calculating second order
sensitivities. The latter quantify the non-linearly coupled effect of
any combination of two parameters. As a completion we used a global sensitivity 
analysis approach from \citet{Chan1997a} to take the variability of parameters into account.
The method we suggest has numerous
uses. A theoretician striving for model reduction may use it for
identifying particularly low sensitivities to detect superfluous parameters. 
An experimenter may use it for identifying particularly high sensitivities to facilitate
determining the parameter value with maximised precision.

We found that, in comparison to Zajac's linear model, Hatze's
non-linear model incorporates a set of parameters to which activation
dynamics is clearly more sensitive. Other than Zajac's model, Hatze's
model can moreover reproduce measured shifts in optimal muscle length
with varied muscle activity. Accordingly, we extracted a specific
parameter set for Hatze's model that combines best with a particular
muscle force-length relation. We also give an outlook on how
sensitivity analysis could be used for optimising parameter sets in
future work.
\end{abstract}
\end{frontmatter}
\clearpage
\section*{List of symbols}
%
%\hfill\\
\begin{table}[h]
\label{Symbol}
\footnotesize
%\renewcommand{\baselinestretch}{1.25}
%\tablesize
\renewcommand{\arraystretch}{0.9}
\begin{tabular}{|c|l|l|}
  \hline
  Symbol & Meaning & Value\\
  \hline\hline
${\ell}_{CE}$		& contractile element (CE) length &  time-depending \\
$\dot{\ell}_{CE}$	& contraction velocity & first time derivative
  of ${\ell}_{CE}$\\
${\ell}_{CEopt}$ 	& optimal CE length & muscle-specific\\
${\ell}_{CErel}$        & relative CE length &
$\ell_{CErel} = \frac{{\ell}_{CE}}{{\ell}_{CEopt}}$ (dimensionless)\\
$F_{max}$               & maximum isometric force of the CE & muscle-specific\\
$\sigma$		& neural muscle stimulation &
 time-depending\,; here: a fixed parameter \\
$q$		        & muscle activity (bound $Ca^{2+}$-concentration) & time-depending \\
$q_0$			& basic activity according to \citet{Hatze2} & $0.005$  \\
$q_H$			& activity according to \citet{Hatze}
& time-length-depending  \\
$q_{H,0}$		& initial condition for Hatze's activation ODE	& mutable  \\
$q_Z$			& activity according to \citet{Zajac}	& time-depending  \\
$q_{Z,0}$		& initial condition for Zajac's activation ODE & mutable  \\
$\tau$			& activation time constant in \citet{Zajac} & here: $\frac{1}{40}\,s$\\
$\tau_{deact}$		& deactivation time constant in \citet{Zajac} &
  here: $\frac{1}{40}\,s$ or $\frac{3}{40}\,s$\\
$\beta$		        & corresponding deactivation boost \citep{Zajac} & $\beta = \tau / \tau_{deact}$ \\
$\nu$			& exponent in Hatze's formulation & 2 or 3 \\
$m$			& activation frequency constant in
  \citet{Hatze} & range: $3.67 \ldots 11.25\,\frac{1}{s}$\,; here: $10\frac{1}{s}$\\
$c$			& maximal
  $Ca^{2+}$-concentration in \citet{Hatze1} & $1.37 \cdot 10^{-4}\,\frac{mol}{l}$\\
$\gamma$	        & representation of free $Ca^{2+}$-concentration \citep{Hatze,Hatze2} 
& time-depending \\
$\rho$			& length dependency of \citet{Hatze1} activation
  dynamics 
  & $\rho ({\ell}_{CErel})= \rho_c \cdot \frac{{\ell}_{\rho}-1}{\frac{{\ell}_{\rho}}{{\ell}_{CErel}}-1}$ \\
$\ell_{\rho}$	        & pole in Hatze's length dependency function & $2.9$ \\
$\rho_0$		& factor in \citet{vanSoest1992a,Hatze} &
$6.62 \cdot 10^{4}\,\frac{l}{mol}$ ($\nu=2$) or $5.27 \cdot 10^{4}\,\frac{l}{mol}$ ($\nu=3$) \\
$\rho_c$		& merging of $\rho_0$ and $c$ &
  $\rho_c=\rho_0 \cdot c$\,; here: $9.10$ ($\nu=2$) or $7.24$ ($\nu=3$)\\
$\Lambda$	        & model parameter set & $\Lambda=\{\lambda_1,\ldots,\lambda_n\}$\\
  \hline
\end{tabular}  
\end{table}
\clearpage
\section{Introduction}
In science, knowledge is gained by an interplay between quantitative
real world measurements of physical, chemical, or biological phenomena
and the development of mathematical models accounting for the
dynamical processes behind. In general, such phenomena are determined
as spatio-temporal patterns of physical measures (state
variables). Modelling consists of distinguishing the surrounding world
from a system that yields the phenomena, and formulating a
mathematical description of the system, a model, that can calculate
its state variables. The calculations depend on model parameters and
often on prescribing measured input variables. By changing parameter
values and analysing the resulting changes in the values of the state
variables, the model may then be used as a predictive tool. This way,
the model's validity can be verified. If the mathematical model
description is moreover derived from first principles, the model can
potentially even explain the phenomena in a causal sense.

Calculating the sensitivities of a model's predicted output, i.e.,
the system's state variables, with respect to model parameters is a
means of eliminating redundancy and indeterminancy from models, and
thus helps to identify valid models. Sensitivity analyses can be
helpful both in model-based experimental approaches and in purely
theoretical work. A modelling theoretician may be looking for
parameters to which all state variables are non-sensitive. Such
parameters may be superfluous. An experimenter may inspect the model
that represents his working hypothesis and analyse which of the
model's state variables are specifically sensitive to a selected
parameter. It would then make sense to measure exactly this state
variable to identify the value of the selected parameter.

In a biomechanical study, \citet{Scovil} applied sensitivity analysis
to examine the dynamics of a mechanical multi-body system (a runner's
skeleton) coupled to muscle activation-contraction dynamics. They
calculated specific sensitivity coefficients in three slightly
different ways. A sensitivity coefficient is the difference quotient
that is calculated from dividing the change in a state variable by the
change in a model parameter value, evaluated in a selected system
state \citep{Lehman1982a}. The corresponding partial derivative may be
simply called ``sensitivity''. A sensitivity function is then the time
evolution of a sensitivity \citep{Lehman1982a}. Thus,
\citet{Lehman1982a} had proposed a more general and unified approach
than \citet{Scovil}, which allows to systematically calculate the
sensitivities of any dynamical system described in terms of ordinary
differential equations. As an example for sensitivity functions,
\citet{Lehman1982a} had applied their proposed method to a
muscle-driven model of saccadic eye movement. By calculating a
percentage change in a state variable value per percentage change in a
parameter value, all sensitivities can be made comprehensively
comparable, even across models. A sensitivity as defined so far is of
first order. Methodically, we aim at introducing a step beyond, namely
at calculating second order sensitivities. These measures are suited
to quantifying how much the sensitivity of a state variable with
respect to (w.r.t.) one model parameter depends on changing another
parameter. This way, the strength of their interdependent influence
on model dynamics can be determined.

In addition to this so-called local sensitivity analysis we furthermore take the variability of the parameters into account aiming for a global sensitivity analysis as presented in \citet{Chan1997a} and \citet{Saltelli2000}. This approach allows to translate sensitivities with respect to parameters into importances of parameters.

In this study, we will apply the sensitivity analysis to models that
predict how the activity of a muscle (its chemical state) changes
when the muscle is stimulated by neural signals (electrical
excitation). Such models are used for simulations of muscles'
contractions coupled to their activation dynamics. Models for coupled
muscular dynamics are often part of neuro-musculo-skeletal models of
biological movement systems. In particular, we want to try and rate
two specific model variants of activation dynamics formulated by
\citet{Zajac} and by \citet{Hatze}. As a first result, we present
an example of a simplified version of the \citet{Zajac} model, in which
sensitivity functions can even be calculated in closed form. Then, we
calculate the sensitivities numerically with respect to all model
parameters in both models, aiming at an increased understanding of the
influence of changes in model parameters on the solutions of the
underlying ordinary differential equations (ODEs). Additionally, we
discuss which of both models may be the more physiological one. The
arguments come from a mixture of three different aspects:
sensitivity analysis, others' experimental findings, and an
additional attempt to best fit different combinations of activation
dynamics and force-length relations of the contractile element (CE) in
a muscle to known data on shifts in optimal CE length with muscle
activity \citep{Kistemaker2005a}.
\section{Two models for muscle activation dynamics} \label{Model}
Macroscopically, a muscle fibre or an assembly of muscle fibres
(muscle belly) is often mapped mathematically by a one-dimensional
massless thread called ``contractile component'' or ``contractile
element'' (CE)
\citep{vanSoest1992a,vanSoest1993b,Cole,GS,Haeufle2014b}. Its
absolute length is ${\ell}_{CE}$ which may be normalised to the
optimal fibre length ${\ell}_{CEopt}$ by
${\ell}_{CErel} = {\ell}_{CE} / {\ell}_{CEopt}$. In macroscopic muscle
models, the CE muscle force is usually modelled as a function of a
\mbox{force-(CE-)}length relation, a \mbox{force-(CE-)}velocity
relation, and \mbox{(CE-)}activity $q$. A common view is that muscle
activity $q$ represents the number of attached cross-bridges within
the muscle, normalised to the maximum number available
($q_{0} \leq q \leq 1$). It can also be considered as
the concentration of bound Ca$^{2+}$-ions in the muscle sarcoplasma
relative to its physiological maximum. The parameter $q_{0}$
represents the minimum activity that is assumed to occur without
any stimulation \citep{Hatze}.

We analyse two different formulations of muscle activation dynamics,
i.e., the time (its symbol: $t$) evolution of muscle activity
$q(t)$. One formulation of muscle activation dynamics was suggested by
\citet{Zajac} which we modified slightly to take $q_{0}$ into account:
\begin{equation}
   \dot{q}_Z = \frac{1}{\tau \cdot (1-q_0)} \cdot \left [ \sigma \cdot
     (1-q_0) - \sigma \cdot (1 - \beta ) \cdot (q_Z-q_0) - \beta \cdot
     (q_Z-q_0) \right ] ,
   \quad q_Z(0)=q_{Z,0} \quad.
 \label{dq_Z} 
\end{equation}
Here, $\sigma$ is meant to represent the (electrical) stimulation of
the muscle, thus, a parameter for controlling muscle dynamics. It
represents the output of the nervous system's dynamics applied to
the muscle which in turn interacts with the skeleton, the body mass
distribution, the external environment, and so with the nervous system
in a feedback loop. Electromyographic (EMG) signals can be seen as a
compound of such neural stimulations collected in a finite volume
(being the input to a number of muscle fibres), over a frequency
range, and coming from a number of (moto-)neurons. The parameter
$\tau$ denotes the activation time constant, and
$\beta = \tau / \tau_{deact}$ is the ratio of activation to
deactivation time constants (deactivation boost).

An alternative formulation of muscle activation dynamics was
introduced by \citet{Hatze}:
\begin{equation} 
   \dot{\gamma} = m \cdot (\sigma - \gamma) \label{Gamma} \quad.
\end{equation}
Here, we divided the original equation from \citet{Hatze} by the parameter
$c = 1.37 \cdot 10^{-4}\,\frac{mol}{l}$
which represents the maximum concentration of free Ca$^{2+}$-ions in the
muscle sarcoplasma. Thus, the values of the corresponding normalised
concentration are $0\leq \gamma \leq 1$.
The activity is then finally calculated by the function
\begin{equation}
    q_H(\gamma,{\ell}_{CErel}) = \frac{q_0 + [ \rho({\ell}_{CErel})
      \cdot \gamma ]^{\nu}}{1 + [ \rho({\ell}_{CErel}) \cdot
      \gamma ]^{\nu}} \label{q_H} \quad,
\end{equation}
and the parameter $c$ is shifted to the accordingly renormalised function
\begin{equation} 
   \rho({\ell}_{CErel}) = \rho_{c} \cdot
   \frac{{\ell}_{\rho} - 1}{\frac{\ell_{\rho}}{{\ell}_{CErel}} - 1}
 \label{rho} \quad,
\end{equation}
with  $\rho_{c} = c \cdot {\rho}_{0}$ and ${\ell}_{\rho} = 2.9$. Two cases
have been suggested by \citet{Hatze2}:
${\rho}_{0} = 6.62 \cdot 10^{4}\,\frac{l}{mol}$
(i.e. $\rho_{c} = 9.10$) for $\nu = 2$
and ${\rho}_{0} = 5.27 \cdot 10^{4}\,\frac{l}{mol}$
(i.e. $\rho_{c} = 7.24$) for  $\nu = 3$ which has been applied in
literature
\citep{vanSoest1992a,Kistemaker2005a,Kistemaker2006a,Kistemaker2007b}.
By substituting equations \eqref{Gamma} and \eqref{q_H} into
$\dot{q_{H}} = \frac{d q_{H}(\gamma,{\ell}_{CErel})}{d \gamma} \cdot \dot{\gamma}$
and resubstituting the inverse of \eqref{q_H} afterwards,
Hatze's formulation of an activation dynamics can be transformed
into a non-linear differential equation directly in terms of the
activity:
\begin{equation}  
    \dot{q}_H = \frac{\nu \cdot m}{1 - q_{0}}\cdot \left [
      \sigma \cdot \rho({\ell}_{CErel}) \cdot (1 - q_{H})^{1+1/\nu}
       \cdot (q_{H} - q_{0})^{1-1/\nu} - (1 - q_{H}) \cdot (q_{H} - q_{0})
 \right ] \label{dq_H} \;, \qquad  \;
\end{equation}
with initial condition $q_H(0)=q_{H,0}$.

The solutions $q_{Z}(t)$ and $q_{H}(t)$ of both formulations of
activation dynamics \eqref{dq_Z} and \eqref{dq_H}, respectively,
can now be directly compared by integrating them with the same initial
condition $q_{Z}(t=0) = q_{H}(t=0)$ using the same stimulation $\sigma$.  
\section{Local first and second order sensitivity of ODE systems regarding their parameters}
\label{Sensitivity}
Let $\Omega \subseteq \R \times \R^M\times \R^N$ and
$f:\, \Omega\rightarrow \R^M$. We then consider a system of ordinary, first
order initial value problems (IVP)
\begin{equation} 
   \dot{Y}=f(t,Y(t,\Lambda),\Lambda) \quad, \qquad Y(0)=Y_0 \quad, \label{ODE}
\end{equation}
where $Y=(y_1(t), y_2(t), \ldots, y_M(t) )$ denotes the vector of
state variables, $f=(f_1,f_2,\ldots,f_M)$ the vector of right-hand
sides of the ODE, and
$\Lambda=\{ \lambda_1, \lambda_2, \ldots , \lambda_N \}$
the set of parameters which the ODE depends on. 
The vector of initial conditions is abbreviated by
\begin{equation}
\label{init}
Y(0)=(y_1(0), y_2(0), \ldots, y_M(0) )=(y_{1,0}, y_{2,0},
\ldots, y_{M,0} )=Y_0 \quad.
\end{equation}
Then, the first order 
solution sensitivity with respect to the parameter set $\Lambda$ is
defined as the matrix
\begin{equation}  
   S(t, \Lambda) = (S_{ik}(t,\Lambda))_{i=1,\ldots, N, k=1,\ldots,M} \quad, \qquad \text{with} \qquad
	S_{ik}(t,\Lambda)=\frac{d}{d \lambda_i} y_k(t,\Lambda) \quad.\label{FO} 
\end{equation} 
For simplicity, we denote $Y=Y(t,\Lambda)$, $f=f(t,Y,\Lambda)$,
$S_{ik}=S_{ik}(t,\Lambda)$ but keep the dependencies in mind. Because
the solution $Y(t)$ might only be gained numerically rather
than in a closed-form expression, we have to apply the well-known
theory of sensitivity analysis as stated in
\citet{Tomovic1972,Dickinson,Lehman1982a,Zivari}. Differentiating
equation \eqref{FO} w.r.t. $t$ and applying the
chain rule yields
\begin{align}
   \frac{d}{dt} S_{ik} &= \frac{d^2}{dt\, d\lambda_i}y_k = \frac{d^2}{d\lambda_i\, dt}y_k
   	= \frac{d}{d\lambda_i}f_k  
   	=\frac{d}{d \lambda_i}Y \cdot  \frac{\partial}{\partial Y}f_k 
   	+ \frac{\partial}{\partial \lambda_i}f_k \quad, \notag
\intertext{with $\frac{\partial}{\partial Y}$ being the gradient of
  state variables. Hence, we obtain the following ODE for the first order
  solution sensitivity}
   \dot{S}_{ik} &= \sum \limits_{l=1}^M  S_{il} \cdot
   \frac{\partial}{\partial y_l}f_k + \frac{\partial}{\partial
     \lambda_i}f_k \quad, \qquad 
   	S_{ik}(0) =\frac{\partial}{\partial \lambda_i} y_{k,0} = 0 \quad, \label{S}
\end{align}
or in short terms
$$\dot{S}=S\cdot J+B \quad, \qquad S(0)=\mathbf{0}_{N\times M} \quad,$$
where $S=S(t)$ is the $N \times M$ sensitivity matrix,
$J=J(t)$ is the $M\times M$ Jacobian matrix with
$J_{kl}=\frac{\partial}{\partial y_l} f_k$, furthermore $B=B(t)$ the
$N\times M$-matrix containing the partial derivatives
$B_{ik}=\frac{\partial}{\partial \lambda_i}f_k$ and
$\mathbf{0}_{N\times M}$ the $N\times M$-matrix consisting of zeros
only.

By analogy, the second order sensitivity of $Y(t)$ with respect to
$\Lambda$ is defined as the following $N \times N \times M$-tensor 
\begin{equation*}
   R(t,\Lambda) = (R_{ijk}(t,\Lambda))_{i,j=1,\ldots N, k=1,\ldots M} \quad,
\end{equation*}
with
\begin{equation}
 \label{S2} 
 R_{ijk}(t,\Lambda) = \frac{d}{d \lambda_i } S_{jk} = \frac{d}{d \lambda_j } S_{ik} 
	= \frac{d^2}{d \lambda_i\, d\lambda_j} y_k = R_{jik}(t,\Lambda) \quad,
\end{equation}
assuming $R_{ijk}=R_{jik}$ for all $k=1,\ldots,M$. That is, we
assume that the prerequisites of Schwarz' Theorem (symmetry of the
second derivatives) are fulfilled throughout.
Again, differentiating w.r.t. $t$ and applying the chain rule
leads to the ODE
\begin{align}
 \label{R} 
 \dot{R}_{ijk}=\sum \limits_{l=1}^M \left(R_{ijl}\frac{\partial}{\partial y_l}f_k 
 + S_{il}\frac{\partial}{\partial \lambda_j}f_k 
 + S_{jl}\frac{\partial}{\partial \lambda_i}f_k \right) 
 +\sum \limits_{l_1=1}^M \sum \limits_{l_2=1}^M 
  S_{il_1}S_{jl_2}\frac{\partial^2}{\partial y_{l_1} \partial y_{l_2}}f_k 
 +\frac{\partial^2}{\partial \lambda_i \partial \lambda_j}f_k \quad,
 \end{align}
with $R_{ijk}(0)=0$. For purposes beyond the aim of this paper, a
condensed notation introducing the concept of tensor (or Kronecker)
products as in \citet{Zivari} may be helpful. For a practical
implementation in \textsc{MatLab} see \citet{Bader}.

Furthermore, if an initial condition $y_{k,0}$ (see \eqref{init}) is
considered as another parameter we can derive a separate sensitivity
differential equation by rewriting equation \eqref{ODE} in its
integral form $$ Y(t) = Y_0 + \int \limits_0^t f(s,Y(s)) \myd s \quad.$$
Differentiating this equation w.r.t. $Y_0$ yields 
$$ S_{Y_0}(t)=\frac{\partial}{\partial Y_0} Y(t)=1+\int
\limits_0^t \frac{\partial }{\partial Y} f \cdot
\frac{\partial}{\partial Y_0} Y(s)  \myd s $$
and differentiating again w.r.t. $t$ results in a homogeneous ODE
for each component $S_{y_{k,0}}(t)$, namely
\begin{equation}
 \dot{S}_{y_{k,0}}(t)=\sum \limits_{l=1}^M \frac{\partial}{\partial
   y_l}f_k \cdot S_{y_{l,0}} \quad, \qquad \text{ with } \qquad S_{y_{k,0}}(0)
   =\frac{\partial}{\partial y_{k,0}} y_{k,0}=1 \quad .\label{y_0} 
 \end{equation}

The parameters of our analysed models are meant to represent
physiological processes, and bear physical dimensions therefore. For
example, $m$ and $\frac{1}{\tau}$ are frequencies measured in [Hz],
whereas $c$ is measured in [mol/l]. Accordingly,
$S_\tau= \frac{d }{d \tau} q_Z$ would be measured in [Hz] and $S_m$ in
[s] (note that our model only consists of $one$ ODE and therefore we
do not need a second index). Normalisation provides a comprehensive
comparison between all sensitivities, even across models. For any
parameter, the value $\lambda_i$ fixed for a specific simulation is a
natural choice. For any state variable, we chose its current value
$y_k(t)$ at each point in time of the corresponding ODE
solution. Hence, we normalise each sensitivity
$S_{ik}=\frac{d y_k}{d\lambda_i}$ by multiplying it with the ratio
$\frac{\lambda_i}{y_k(t)}$ to get the relative
sensitivity
\begin{equation}  
\tilde{S}_{ik} = S_{ik} \cdot \frac{\lambda_i}{y_k} \quad.\label{Snorm} 
\end{equation}
A relative sensitivity $\tilde{S}_{ik}$ thus quantifies the percentage
change in the $k$-th state variable value per percentage change in the
$i$-th parameter value.
This applies accordingly to the second order sensitivity
\begin{equation}  
\tilde{R}_{ijk} = R_{ijk} \cdot \frac{\lambda_i \cdot \lambda_j}{y_k}
\quad. \label{S2norm} 
\end{equation}
It can be shown that this method is valid and mathematically
equivalent to another common method in which the whole model is
non-dimensionalised a priori \citep{Scherzer}. A non-normalised model
formulation has the additional advantage of usually allowing a more
immediate appreciation of and transparent access for experimenters. In
the remainder of this  manuscript, we always present and discuss
relative sensitivity values normalised that way.

In our model, the specific case $M=1$ applies, so 
equations \eqref{S} and \eqref{R} simplify to the case $k=1$ (no summation).
\section{Variance-based global sensitivity analysis}
\label{Variance}
The previous presented differential sensitivity analysis is called a local method because 
it does not take the physiological range of parameter values into account. If we imagine 
the parameter space as a $N$-dimensional cuboid 
$\mathcal{C}=[\lambda_1^-;\lambda_1^+] \times \ldots \times [\lambda_N^-;\lambda_N^+]$, 
where $\lambda_i^-,\lambda_i^+$ are the minimal and maximal parameter value, we 
can only fix a certain point $\hat{\Lambda}=(\hat{\lambda}_1,\ldots, \hat{\lambda}_N)\in \mathcal{C}$ 
and calculate the local gradient of the solution w.r.t. $\hat{\Lambda}$. By changing only 
one parameter at once the investigated star-shaped area lies within a ball around 
$\hat{\Lambda}$ whose volume vanishes in comparison to $\mathcal{C}$ for an 
increasing number of parameters as shown in \citet{Saltelli2010a}.

For taking the range of parameter values into account, \citet{Saltelli2000} gave a detailed 
elaboration of so-called global methods. The main idea behind most global methods is 
to include a statistical component to scan the whole parameter space $\mathcal{C}$ and 
combine the percentage change of the state variable per percentage change of the 
parameters with the variability of the parameters themselves. For an overview of the 
numerous methods like ANOVA, FAST, Regression or Sobol' Indexing we refer the reader 
to \citet{Saltelli2000} and \citet{Frey}.

In this paper we want to sketch the main idea of the variance-based sensitivity analysis 
approach presented in \citet{Chan1997a} based on Sobol' indexing. We chose this method 
because of its transparency and low computational cost. The aim of this method is to 
calculate two measurands of sensitivity w.r.t. parameter $\lambda_i$: the variance based 
sensitivity function denoted by $VBS_i(t)$ and the total sensitivity index function denoted by 
$TSI_i(t)$. The $VBS$ functions give a normalised first order sensitivity quite similar to $\tilde{S}$ 
from the previous section but include the parameter range. The $TSI$ functions, however, 
even include higher order sensitivities and give a measurand for interactions of parameter influences. 

A receipt for calculating $VBS$ and $TSI$ can be given as follows. First of all set boundaries 
for all model parameters, either by model assumptions or literature reference, thus by 
resulting in $\mathcal{C}$. Secondly  generate two sets of $n$ sample points 
$\hat{\Lambda}_{1,j},\hat{\Lambda}_{2,j} \in \mathcal{C}, \;j=1,..,n$ w.r.t to the underlying 
probability distribution of each parameter, in our case the uniform distribution. Thirdly 
calculate $2nN$ sets of new sample points 
$\hat{\Lambda}_{1,j}^i,\hat{\Lambda}_{1,j} ^{\sim i}, \; j=1,..,n, \; i=1,..,N$ where 
$\hat{\Lambda}_{1,j}^i$ consists of all sample points in $\hat{\Lambda}_{1,j}$ with the 
$i$-th component of $\hat{\Lambda}_{2,j}$. Consequently $\hat{\Lambda}_{1,j} ^{\sim i}$ 
consists of the $i$-th component of $\hat{\Lambda}_{1,j}$ and every other component of 
$\hat{\Lambda}_{2,j}$. Fourthly evaluate the model from Eqn.~\eqref{ODE} on every of the 
$2n (N+1)$ sample points $\hat{\Lambda}_{1,j},\hat{\Lambda}_{2,j} ,\hat{\Lambda}_{1,j}^i, 
\hat{\Lambda}_{1,j} ^{\sim i}$ resulting in a family of solutions. %$\mathcal{F}[Y](t)$
For this family perform the following calculations:
\begin{enumerate}
\item Compute the variance of the family of all $2n (N+1)$ solutions as a function of time, 
namely $V(t)$. This variance function indicates the general model output variety throughout
 the whole parameter range.
\item Compute the variance of the family of all $nN+1$ solutions resulting from an evaluation 
of the model at $\hat{\Lambda}_{1,j} $ and $\hat{\Lambda}_{1,j}^i$. Again the variance is a 
function of time, namely $V_i(t)$, that indicates the model output variety if only one parameter 
value is changed.
\item Compute the variance of the family of all $nN+1$ solutions resulting from an evaluation 
of the model at $\hat{\Lambda}_{1,j}$ and $\hat{\Lambda}_{1,j}^{\sim i}$, namely 
$V_{\sim i}(t)$, that indicates the model output variety if only one parameter value fixed.
\end{enumerate}
Note that in \citet{Chan1997a} the coputations are done via an approximation by Monte-Carlo integrals. The $VBS$ and $TSI$ can finally be calculated as 
\begin{equation} \label{VBS}
VBS_i(t)=\frac{V_i(t)}{V(t)} , \qquad TSI_i= 1-\frac{V_{\sim i}(t)}{V(t)}
\end{equation}
As a consequence of the normalisation we can give additional properties of $VBS$ and 
$TSI$ that can be comprehended in \citep[Fig. 1]{Chan1997a}:
\begin{equation} \label{VBS_props}
\sum\limits_{i=1}^N  VBS_i(t)\leq 1 , \qquad \sum\limits_{i=1}^N TSI_i(t)\geq 1
\end{equation}
In other words $VBS_i(t)$ gives the normalised first order sensitivity function of the 
solution w.r.t. $\lambda_i$ in relation to the model output range. Accordingly $TSI_i(t)$ 
gives a relative influence of parameter $\lambda_i$ on the model output regarding all 
interactions between other parameters. \citet{Chan1997a} suggested relating the $TSI_i(t)$ 
value to the ``importance" of $\lambda_i$.

\section{An analytical example for local sensitivity analysis including a link between Zajac's and Hatze's
  formulations}
\label{Example}
By further simplifying Zajac's formulation of an activation dynamics
\eqref{dq_Z} through assuming an deactivation boost $\beta=1$
(activation and deactivation time constants are equal)
and a basic activity $q_0=0$, we obtain a linear ODE for this specific
case $q_Z^{sp}$ which is equivalent to Hatze's equation \eqref{Gamma}
modelling the time evolution of the free Ca$^{2+}$-ion concentration:
\begin{equation}
\dot{ q}_Z^{sp}=\frac{1}{\tau}(\sigma -q_Z^{sp}) \quad, \qquad q_Z^{sp}(0)=q_{Z,0} 
\quad .\label{special}
\end{equation}
By analysing this specific case, we aim at making the above described
sensitivity analysis method more transparent for the reader. Solving
equation \eqref{special} yields
\begin{equation}
q_Z^{sp}(t)= \sigma \cdot (1 - e^{-t/ \tau}) +q_{Z,0} \cdot e^{-t/\tau}  \label{specialsolution}
\end{equation}
depending on just two parameters $\sigma$ (stimulation: control
parameter) and $\tau$ (time constant of activation: internal
parameter) in addition to the initial value $y_{0} = q_{Z,0}$.
The solution $q_{Z}(t)$ equals the $\sigma$ value after about $\tau$.

\begin{figure}[h!]
\includegraphics[width=10cm]{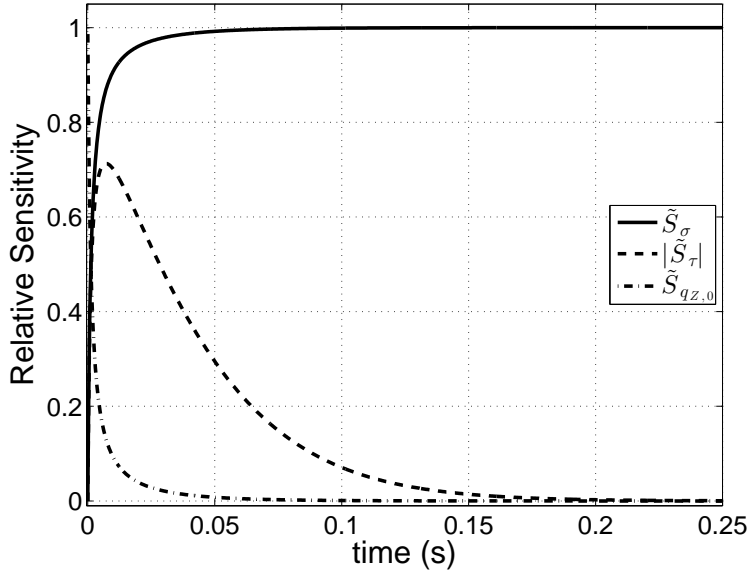}
\caption{Relative sensitivities $\tilde{S}_{i}$ w.r.t. the three parameters
  in the simplified formulation \eqref{special} of Zajac's activation
  dynamics \eqref{dq_Z}. Parameters: stimulation
  $\sigma$ (see \eqref{S_special_sigma}: solid line), activation
  time-constant $\tau$ (see \eqref{S_special_tau}:
  dashed line), and initial activation $q_{Z,0}$ (see
  \eqref{S_special_qZ0}: dash-dotted line). Note that
  $\tilde{S}_{\tau}$ is negative, but for reasons of comparability we
  have plotted its absolute value. Parameter values are
  $\sigma=1$, $\tau = \frac{1}{40}\,s = 0.025\,s$, and
  $q_{Z,0}=0.05$. Because the ODE
  \eqref{special} for $q_Z^{sp}$ is equivalent to Hatze's ODE
  \eqref{Gamma} for the free Ca$^{2+}$-ion concentration, $\gamma$, we
  can identify the sensitivity of $\frac{1}{\tau}$ with that of $m$.
  \label{fig_analytic}}
\end{figure}

Although already knowing the solution \eqref{specialsolution}
explicitly in terms of time and both parameters, we still apply,
because of transparency, the more generally applicable, implicit
method \eqref{S},\eqref{y_0} to determine the derivatives of the
solution w.r.t. the parameters (the sensitivities). For that, we
calculate the gradient of the right hand side $f(q_Z^{sp}, \sigma,
\tau)$ of the ODE
\eqref{special}

\begin{equation*}
   \frac{\partial}{\partial q_Z^{sp}}f = -\frac{1}{\tau}    \quad , \quad 
   \frac{\partial}{\partial \sigma}f =     \frac{1}{\tau}  \quad,
   \text{and} \quad
   \frac{\partial}{\partial \tau}f =- \frac{\sigma-q_Z^{sp}}{\tau^2}
   = \frac{q_{Z,0}-\sigma}{\tau^2} e^{-t/ \tau}
\end{equation*}
and insert these partial derivatives into equations
\eqref{S} and \eqref{y_0}. Solving the respective three ODEs for the
three parameters ($\sigma, \tau$, $q_{Z,0}$) and
normalising them according to \eqref{Snorm} gives the relative
sensitivities of $q_Z^{sp}$ w.r.t. $\sigma$, $\tau$, and $q_{Z,0}$ as
functions of time (see Fig.\,\ref{fig_analytic}):
\begin{equation}
\tilde{S}_\sigma(t)=(1-e^{-t/\tau})\cdot \frac{\sigma}{q_Z^{sp}(t)} 
=\frac{ \sigma \cdot ( e^{t/\tau}-1) }{ \sigma \cdot ( e^{t/\tau}-1) +q_{Z,0} }
\quad ,\label{S_special_sigma}
\end{equation}
\begin{equation}
\tilde{S}_\tau (t)= \left(\frac{(q_{Z,0}-\sigma) \cdot
  t }{\tau^2}e^{-t/\tau}\right) \cdot \frac{\tau}{q_Z^{sp}(t)}=
\frac{t\cdot (q_{Z,0}-\sigma )}{\tau \cdot [ \sigma \cdot ( e^{t/\tau}-1) +q_{Z,0} ]} 
\quad \text{, and}  \label{S_special_tau}
\end{equation}
\begin{equation}
\tilde{S}_{q_{Z,0}}(t) = e^{-t/\tau} \cdot \frac{q_{Z,0} }{q_Z^{sp}(t)}=\frac{q_{Z,0}}{  
\sigma \cdot ( e^{t/\tau}-1) +q_{Z,0} } \quad. \label{S_special_qZ0}
\end{equation}

A straight forward result is that the time constant $\tau$ has its
maximum effect on the solution (Fig.\,\ref{fig_analytic}: see
$\tilde{S}_{\tau}(t)$) at time $t = \tau$. In case of a step in
stimulation, the sensitivity $\tilde{S}_{\tau}(t)$ vanishes in the
initial situation and exponentially approaches zero again after a few
further times the typical period $\tau$. Note that
$\tilde{S}_{\tau}(t)$ is negative which means that an increase in
$\tau$ decelerates activation. That is, for a fixed initial value
$q_{Z,0}$, the solution value $q_{Z}(t)$ decreases at a given point in
time if $\tau$ is increased. After a step in stimulation $\sigma$, the
time in which the solution $q_{Z}(t)$ bears some memory of its initial
value $q_{Z,0}$ is equal to the period of being non-sensitive to any
further step in $\sigma$ (compare $\tilde{S}_{q_{Z,0}}(t)$ to
$\tilde{S}_{\sigma}(t)$ and \eqref{S_special_sigma} to
\eqref{S_special_qZ0}). After about $\tau/2$ the sensitivity
$\tilde{S}_{q_{Z,0}}(t)$ has already fallen to about $0.1$ and
$\tilde{S}_{\sigma}(t)$ to about $0.9$ accordingly.
\section{The numerical approach and results} \label{Num}
Typically, biological dynamics are represented by non-linear ODEs. So the
linear ODE used for describing activation dynamics in the
\citet{Zajac} case \eqref{dq_Z} is more of an exception. For example,
a closed-form solution can be given. \eqref{specialsolution} is
an example as shown in the previous section for
the reduced case of non-boosted deactivation \eqref{special}.

In general, however, non-linear ODEs used in biomechanical modelling,
as the \citet{Hatze} case \eqref{dq_H} for describing activation
dynamics, can only be solved numerically. It is understood that any
explicit formulation of a model in terms of ODEs allows to provide the
partial derivatives of their right hand sides $f$ w.r.t. the model
parameters in a closed form.
Fortunately, this is exactly what is required as part of the
sensitivity analysis approach presented in section
\ref{Sensitivity}, in particular in equation \eqref{S}.

As an application for applying this approach, we will now present a
comparison of both formulations of activation dynamics. The example
indicates that the approach may be of general value
because it is common practice in biomechanical modelling to (i)
formulate the ODEs in closed form and (ii) integrate the ODEs
numerically. Adding further sensitivity ODEs for model parameters
is then an inexpensive enhancement of the procedure used to solve
the problem anyway.

For the two different activation dynamics \citep{Zajac} and
\citep{Hatze}, the parameter sets $\Lambda_Z$ and $\Lambda_H$,
respectively, consist of
\begin{align} 
   \Lambda_Z &= \{\,q_{Z,0},\,\sigma,\,q_0,\,\tau,\,\beta\,\} \quad, \label{lam_Z}\\
   \Lambda_H &= \{\,q_{H,0},\,\sigma,\,q_0, m,\,\rho_c,\,\nu,\,\ell_{\rho},\,{\ell}_{CErel}\,\} 
   \quad, \label{lam_H} 
    \end{align}
including the initial conditions.
The numerical solutions for these ODEs were
computed within the MATLAB environment (The MathWorks,
Natick, USA; version R2013b) using
the pre-implemented numerical solver $ode45$ which is a
Runge-Kutta algorithm of order 5 (for details see \citep{DGL}).
\subsection{Results for Zajac's activation dynamics: sensitivity
  functions}
\label{Zajac_result}
We simulated activation dynamics for the parameter set $\Lambda_Z$
\eqref{lam_Z} leaving two of the values constant ($q_0=0.005$,
$\tau=\frac{1}{40}\,s$) and varying the other three (initial condition
$q_{Z,0}$, stimulation $\sigma$, and deactivation boost $\beta$). The
time courses of the relative sensitivities $\tilde{S}_{i}(t)$
w.r.t. all parameters $\lambda_{i} \in \Lambda_Z$
are plotted in Fig.\,\ref{fig_Zajac_result}.
In the left column of Fig.\,\ref{fig_Zajac_result}, $\beta=1$ is used,
in the right column $\beta=1/3$. Pairs of the parameter values
$q_0=0.005 \leq q_{Z,0} \leq 0.5$ and $0.01 \leq \sigma \leq 1$ are
specified in the legend of Fig.\,\ref{fig_Zajac_result}, with
increasing values of both parameters from top to bottom.
\subsubsection*{Relative sensitivity $\tilde{S}_{q_0}$}
Solutions are non-sensitive to the $q_0$ choice except if both initial
activity and stimulation (also approximating the final activity if
$\beta=1$ and $\sigma >> q_0$) are very low nearby $q_0$ itself.
\subsubsection*{Relative sensitivity $\tilde{S}_{q_{Z,0}}$}
The memory (influence on solution) of the initial value is lost
after about $2 \tau$, almost independently of all other
parameters. This loss in memory is obviously slower than in the
extreme case $q_{Z,0}=0$ (initial value) and $\sigma=1$ (for $\beta=1$
and $q_0=0$ exactly the final value; see section \ref{Example} and
Fig.\,\ref{fig_analytic}). In that extreme case, the influence
(relative sensitivity) of the lowest possible initial value
($q_{Z,0}=0$) on the most rapidly increasing solution (maximum
possible final value: $\sigma=1$) is lost earlier.
\subsubsection*{Relative sensitivity $\tilde{S}_{\tau}$}
The influence of the time constant $\tau$ on the solution is
reduced with decreasing difference between initial and final activity
values (compare maximum $\tilde{S}_{\tau}$ values in
Figs.\,\ref{fig_analytic} and \ref{fig_Zajac_result}) and, no matter
the $\beta$ value, with compoundly raised levels of initial activity
$q_{Z,0}$ and $\sigma$, the latter determining the final activity
value if $\beta=1$. When deactivation is slower than activation
($\beta<1$: right column in Fig.\,\ref{fig_Zajac_result})
$\tilde{S}_{\tau}$ is higher than in the case $\beta=1$, both in its
maximum amplitude and for longer times after the step in stimulation,
especially at low activity levels (upper rows in
Fig.\,\ref{fig_Zajac_result}).
\subsubsection*{Relative sensitivity $\tilde{S}_{\sigma}$}
Across all parameters, the solution is in general most sensitive to
$\sigma$. However, the influence of the deactivation boost parameter
$\beta$ is usually comparable. In some situations, this holds also
for the activation time constant $\tau$ (see below). For $\beta=1$
(Fig.\,\ref{fig_Zajac_result}, left), the solution becomes a
little less sensitive to $\sigma$ with decreasing activity level
($\tilde{S}_{\sigma}<1$), which reflects that the final
solution value is not determined by $\sigma$ alone but by $q_0>0$
and $\beta \neq 1$ as much. If deactivation is much slower than
activation ($\beta=\frac{1}{3}<1$: Fig.\,\ref{fig_Zajac_result},
right), we find the opposite to the case $\beta=1$: $\sigma$
determines the solution the less the more the activity level
rises. Additionally, stimulation $\sigma$ somehow competes with both
deactivation boost $\beta$ and time constant $\tau$ (see further
below). Using the term ``compete'' is meant to illustrate the idea
that any single parameter should have in a sense an individual
interest in influencing the dynamics as much as
possible in order not to be considered superfluous.
\subsubsection*{Relative sensitivity $\tilde{S}_{\beta}$}
Sensitivity w.r.t $\beta$ generally decreases with increasing activity
$q_{Z,0}$ and stimulation $\sigma$ levels, and vanishes at maximum
stimulation $\sigma=1$.
\subsubsection*{Relative sensitivities $\tilde{S}_{\sigma}$,
  $\tilde{S}_{\beta}$, $\tilde{S}_{\tau}$}
\begin{figure}[!ht]
$\begin{matrix} 
\includegraphics[width=6cm]{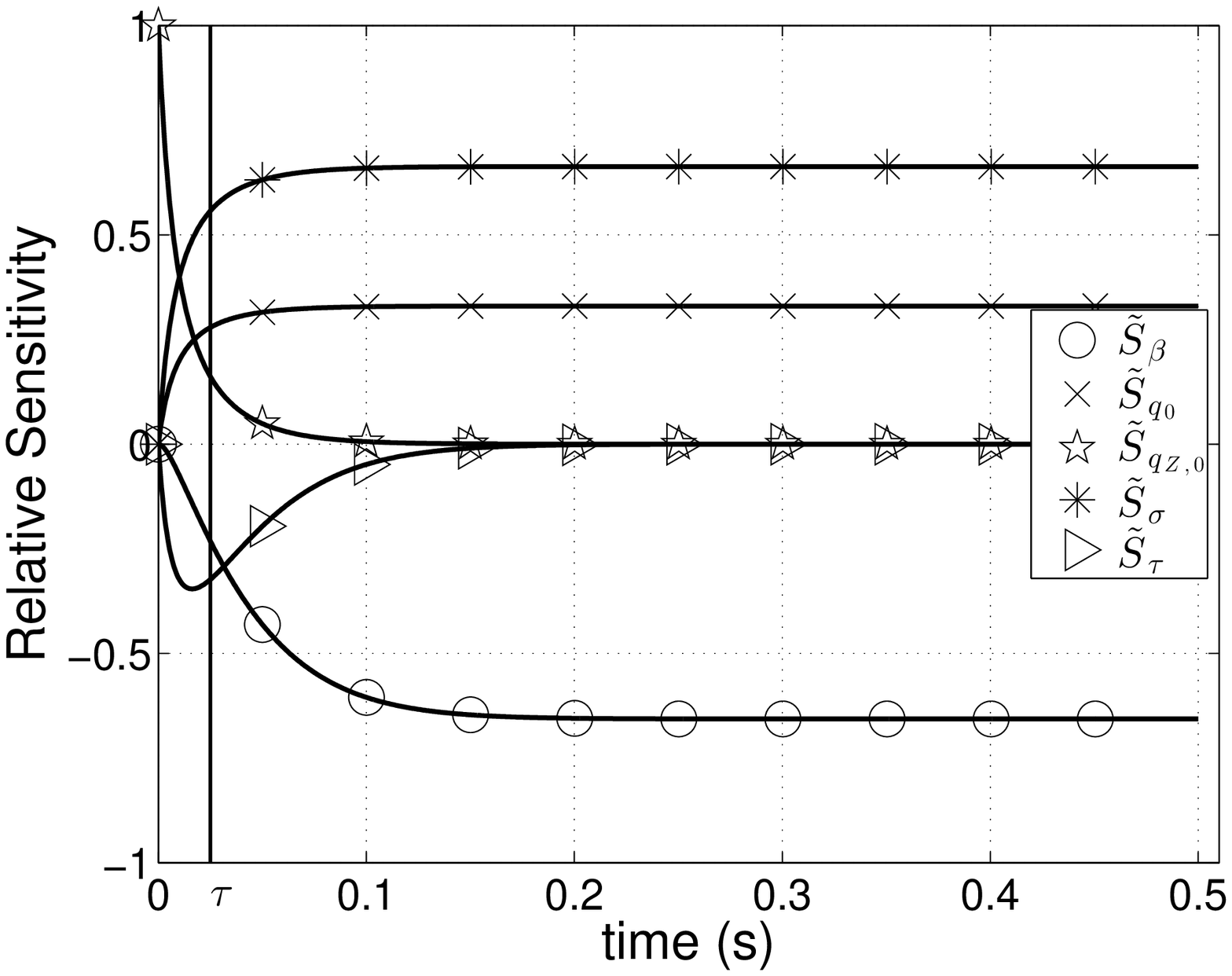} & \includegraphics[width=6cm]{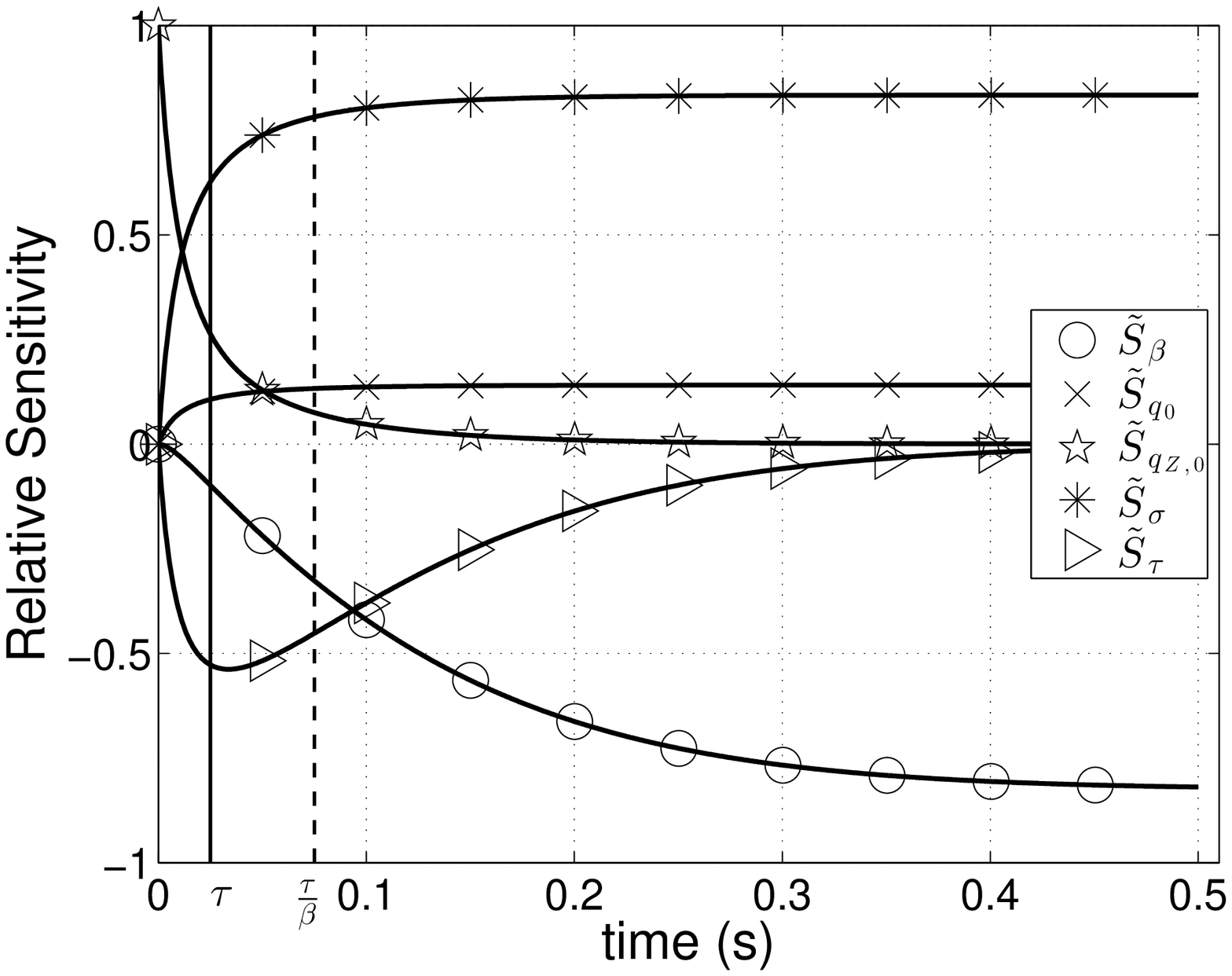} \\
\includegraphics[width=6cm]{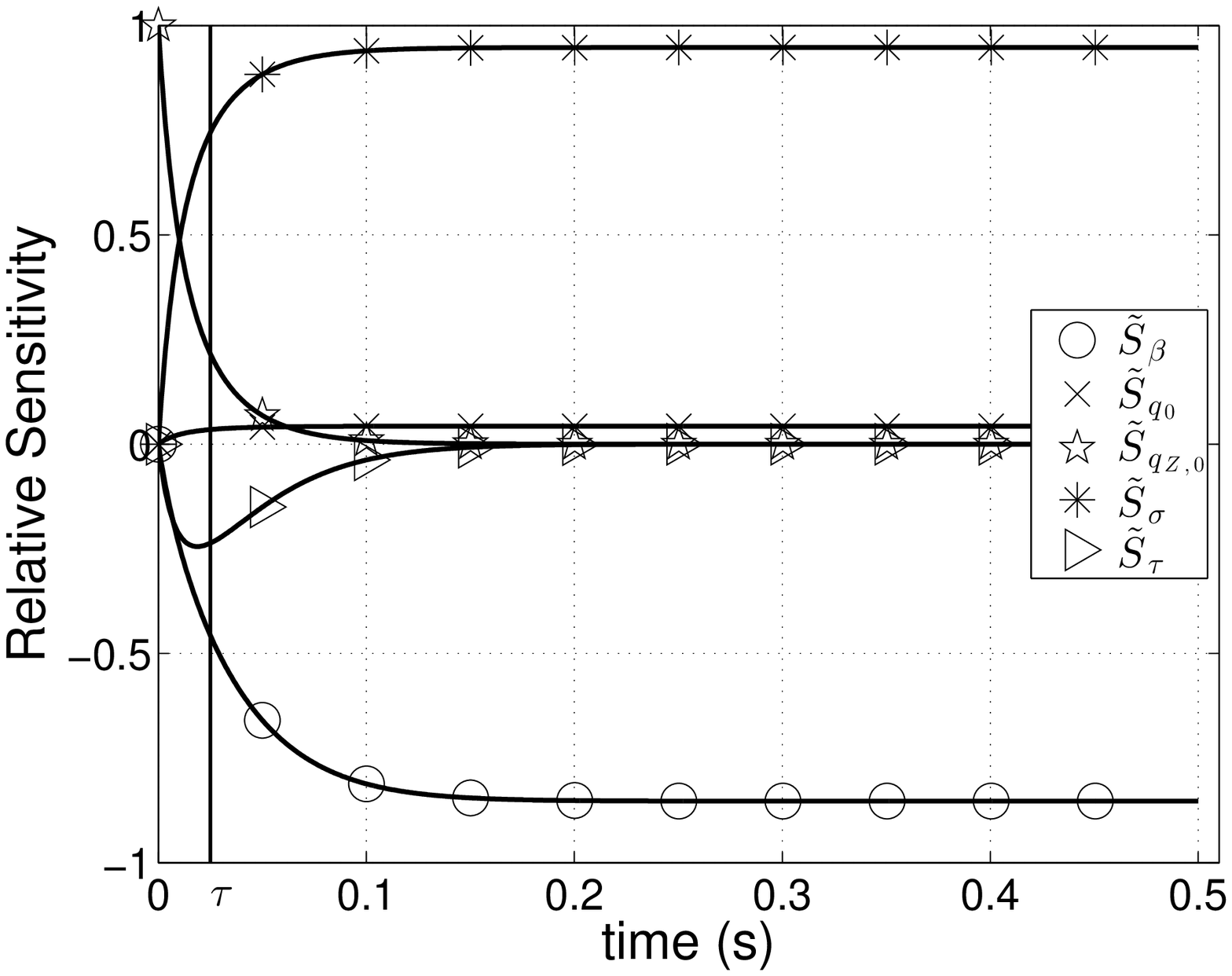} &\includegraphics[width=6cm]{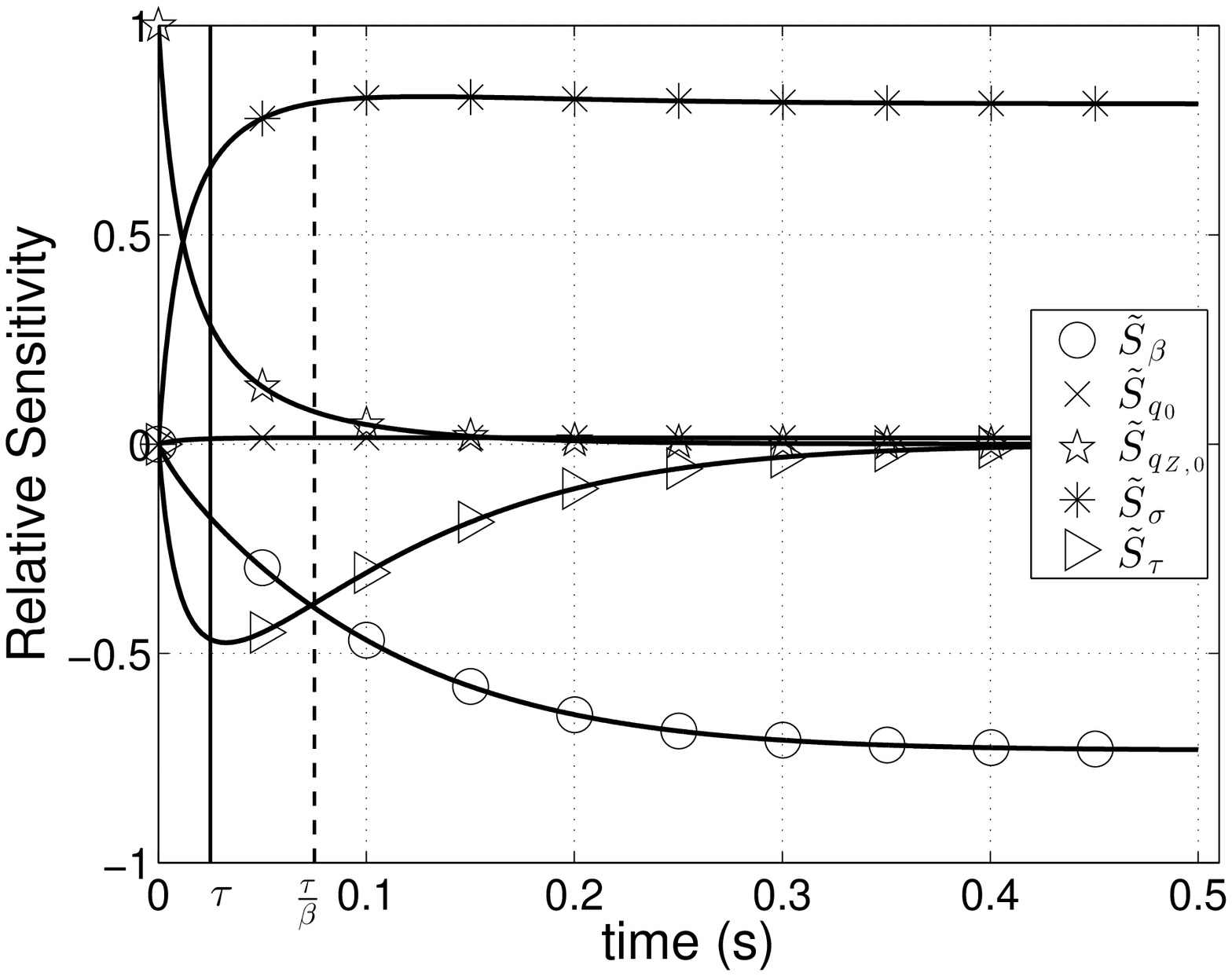} \\
\includegraphics[width=6cm]{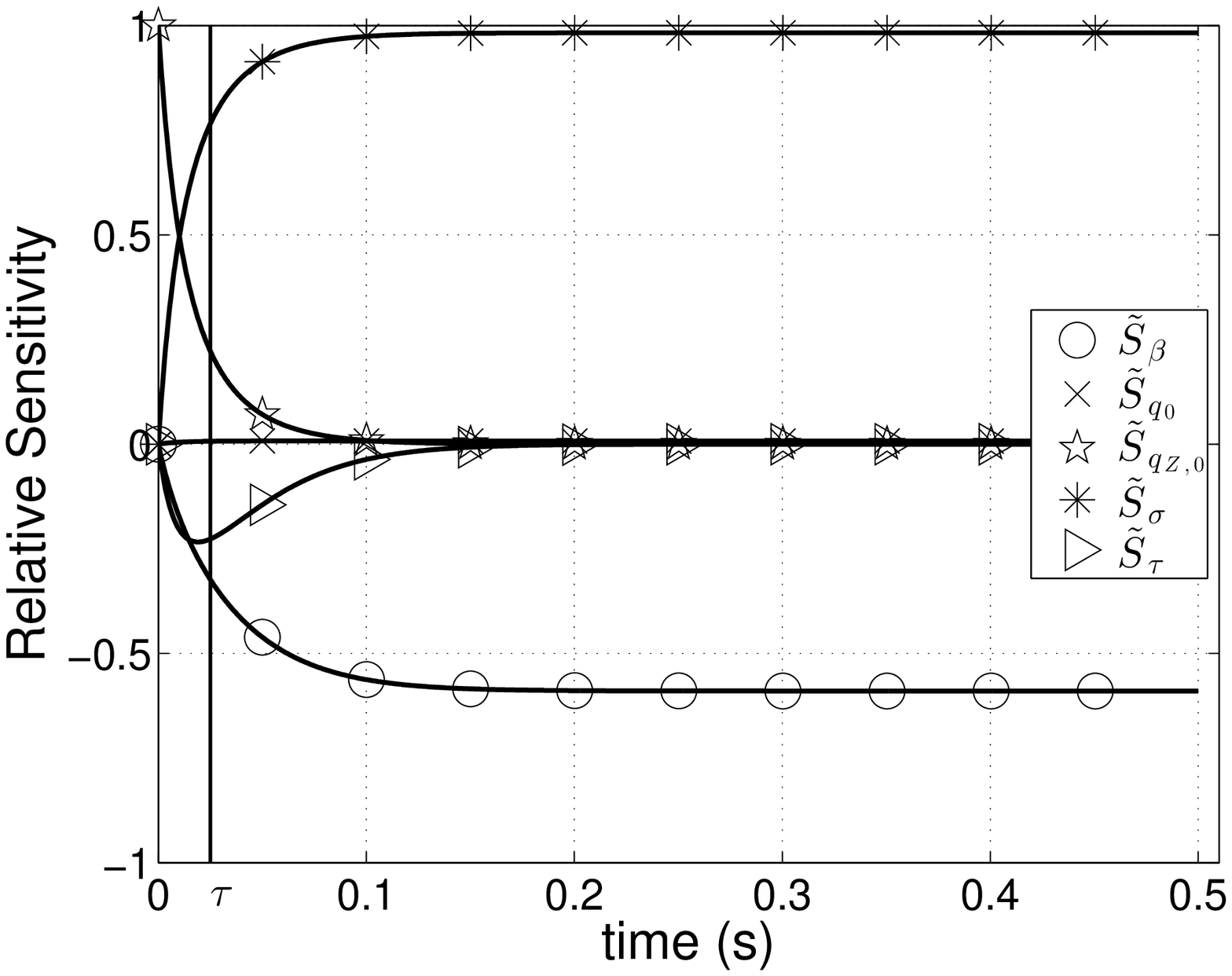} & \includegraphics[width=6cm]{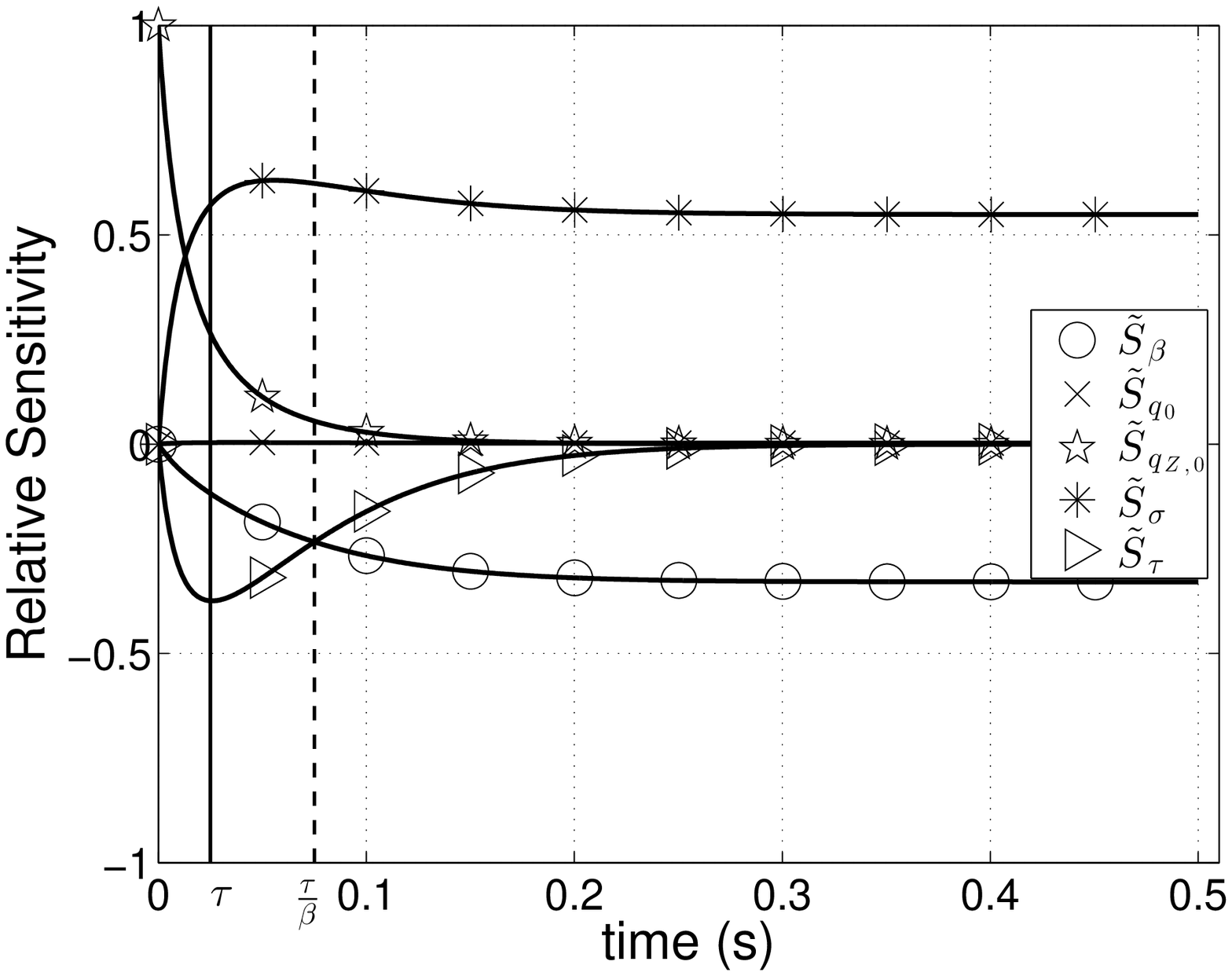} \\
\includegraphics[width=6cm]{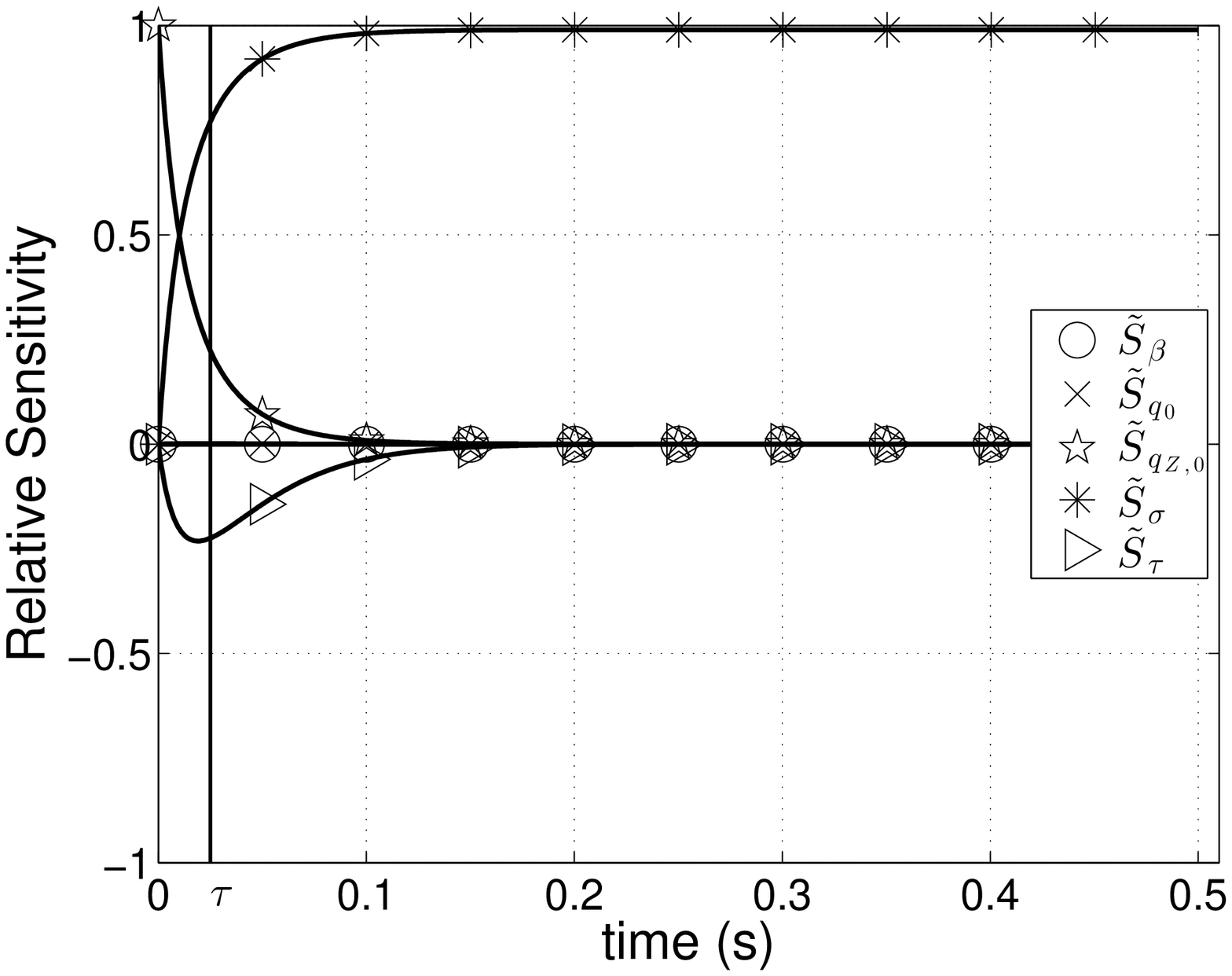} & \includegraphics[width=6cm]{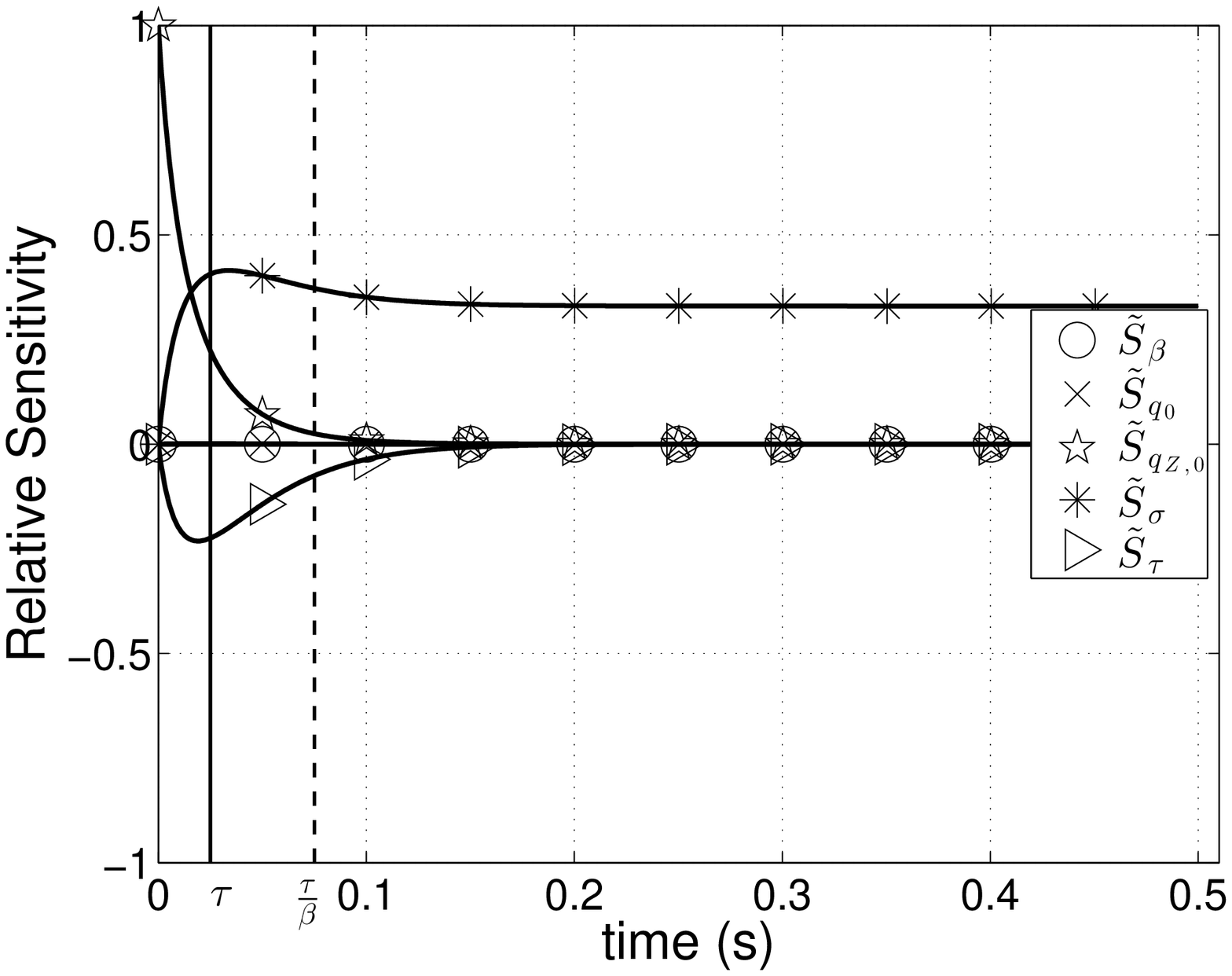} 
\end{matrix}$
\caption{Relative sensitivities $\tilde{S}_{i}$ w.r.t. all parameters
  $\lambda_i$ (set $\Lambda_Z$ \eqref{lam_Z}) in Zajac's activation
  dynamics \eqref{dq_Z}. Parameter values varied from top (i) to
  bottom (iv) row:
\mbox{(i) $q_{Z,0}=q_0=0.005,\,\sigma = 0.01$}, 
\mbox{(ii) $q_{Z,0}=0.05,\,\sigma = 0.1$},
\mbox{(iii) $q_{Z,0}=0.2,\,\sigma = 0.4$},
\mbox{(iv) $q_{Z,0}=0.5,\,\sigma = 1$}; left column: $\beta=1$, right
column:  $\beta=1/3$. \label{fig_Zajac_result}}
\end{figure}

At submaximal stimulation levels $\sigma<1$, the final solution value
is determined to almost the same degree by stimulation $\sigma$ and
deactivation boost $\beta$, yet, with opposite tendencies
($\tilde{S}_{\sigma}>0$, $\tilde{S}_{\beta}<0$). Both parameters
compete, in the above explained meaning, for their impact on the final
solution value. Only at maximum stimulation $\sigma=1$ (lowest row in
Fig.\,\ref{fig_Zajac_result}), this
parameter competition is resolved in favour of $\sigma$. In this
specific case, $\beta$ does not influence the solution at all.
For $\beta=1$ the competition about influencing the solution is
intermittently but only slightly biased by $\tau$: sensitivity
$\tilde{S}_{\tau}$ peaks at comparably low magnitude around
$t=\tau$. This $\tau$ influence comes likewise intermittently at the
cost of $\beta$ influence: the absolute value of
$\tilde{S}_{\beta}$ rises a little slower than
$\tilde{S}_{\sigma}$. In the case $\beta<1$, this competition
becomes even more differentiated and spreaded out in time. Again at
submaximal stimulation and activity levels, the absolute value of
$\tilde{S}_{\tau}$ is lower than that of $\tilde{S}_{\sigma}$ but
higher than that of $\tilde{S}_{\beta}$, making all three parameters
$\sigma$, $\beta$, and $\tau$ compete to comparable degrees for an
impact on the solution until about $t=4 \tau$. Also,
$\tilde{S}_{\tau}$ does not vanish before about $t=10 \tau$.
\clearpage
\subsection{Results for Hatze's activation dynamics: sensitivity
  functions}
\label{Hatze_result}
We also simulated activation dynamics for the parameter set
$\Lambda_H$ \eqref{lam_H}, leaving now four of the values constant
($q_0=0.005$, $m=10\,\frac{1}{s}$, $\ell_{\rho}=2.9$, $\ell_{CErel}=1$)
and again varying three others (initial condition $q_{Z,0}$, stimulation
$\sigma$, and non-linearity $\nu$), keeping in mind that the eighth
parameter ($\rho_c$) is assumed to depend on $\nu$. Again, time
courses of the relative sensitivities $\tilde{S}_{i}(t)$ w.r.t. all
parameters $\lambda_{i}$ (set $\Lambda_H$) are plotted (see
Fig.\,\ref{fig_Hatze_result}). In the left column of
Fig.\,\ref{fig_Hatze_result}, $\nu=2,\,\rho_c=9.10$ is used, in the
right column $\nu=3,\,\rho_c=7.24$. Here again, the same pairs of the
parameter values
($q_0=0.005 \leq q_{Z,0} \leq 0.5$ and $0.01 \leq \sigma \leq 1$,
increasing from top to bottom; see legend of
Fig.\,\ref{fig_Hatze_result}) are used as in the previous
section \ref{Zajac_result} (Fig.\,\ref{fig_Zajac_result}).

Hatze's activation dynamics \eqref{dq_H} are non-linear unlike
Zajac's activation dynamics \eqref{dq_Z}. This non-linearity
manifests particularly in a changeful influence of the parameter $\nu$.
Additionally, the parameter $m$ is just roughly comparable to the
inverse of the exponential time constant $\tau$ in Zajac's linear
activation dynamics.

\subsubsection*{Relative sensitivity $\tilde{S}_{m}$}
In Zajac's linear differential equation \eqref{dq_Z}, $\tau$
establishes a distinct time scale independent of all other parameters.
The parameter $m$ in Hatze's activation dynamics \eqref{dq_H} is just
formally equivalent to the reciprocal of $\tau$: the
sensitivity $\tilde{S}_{m}$ does not peak stringently at
$t=1/m=0.1\,s$ but rather diffusely between about $0.05\,s$ and
$0.1\,s$ in both cases $\nu=2$ and $\nu=3$. At first sight, this is
not a surprise because the scaling factor in Hatze's dynamics is
$\nu \cdot m$ rather than just $m$. However, $\nu \cdot m$ does
neither fix an invariant time scale for Hatze's non-linear
differential equation. This fact becomes particularly
prominent at extremely low activity levels for $\nu=2$
(Fig.\,\ref{fig_Hatze_result}, left, top row) and up to moderately
submaximal activity levels for $\nu=3$
(Fig.\,\ref{fig_Hatze_result}, right, top two rows). Here,
$\tilde{S}_{m}$ is negative which means that increasing the parameter
$m$ results in less steeply increasing activity. This observation is
counter-intuitive to identifying $m$ with a reciprocal of a time
constant like $\tau$. Moreover, as might be expected from the product
$\nu \cdot m$, the exponent $\nu$ does not linearly scale the time
behaviour because $\tilde{S}_{m}$ peaks do not occur systematically
earlier in the $\nu=3$ case as compared to $\nu=2$.
\subsubsection*{Relative sensitivity $\tilde{S}_{q_{H,0}}$}
Losing the memory of the initial condition confirms the above
analysis of time behaviour based on $\tilde{S}_{m}$. At high
activity levels (Fig.\,\ref{fig_Hatze_result}, bottom row),
Hatze's activation dynamics loses memory after practically identical
time horizons, no matter the $\nu$ value, seemingly even slower for
higher $\nu$ at intermediate levels (Fig.\,\ref{fig_Hatze_result}, two
middle rows), and clearly faster at very low
levels (Fig.\,\ref{fig_Hatze_result}, top row). Still, the parameter
$m$ does roughly determine the time horizon in which the
memory of the initial condition $q_{H,0}$ is lost and the
influence of all other parameters is continuously
switched on from zero influence at $t=0$. 
\subsubsection*{Relative sensitivity $\tilde{S}_{q_0}$}
As in Zajac's dynamics the solution is generally only sensitive to
$q_0$ at very low stimulation levels $\sigma \approx q_0$
(Fig.\,\ref{fig_Hatze_result}, top row). Only at
such levels, the $\nu=3$ case shows, however, the peculiarity that the
solution becomes strikingly insensitive to any other parameter than
$q_0$ itself (and $q_{H,0}$). The time evolution of the solution is
more or less determined by just this minimum ($q_0$) and initial
($q_{H,0}$) activities, and $m$ determining the approximate switching
time horizon between both. In particular, the $\ell_{CE}$ dependency,
constituting a crucial property of Hatze's
activation dynamics, is practically suppressed for $\nu=3$ at very low
activities and stimulations. For $\nu=2$, in contrast,
$\tilde{S}_{\ell_{CErel}}$ remains on a low but still significant level
of about a fourth of the three dominating quantities
$\tilde{S}_{q_0}$, $\tilde{S}_{q_{H,0}}$, and $\tilde{S}_{\nu}$.
\subsubsection*{Relative sensitivity $\tilde{S}_{\nu}$}
\begin{figure}[!ht]
$\begin{matrix} 
\includegraphics[width=6cm]{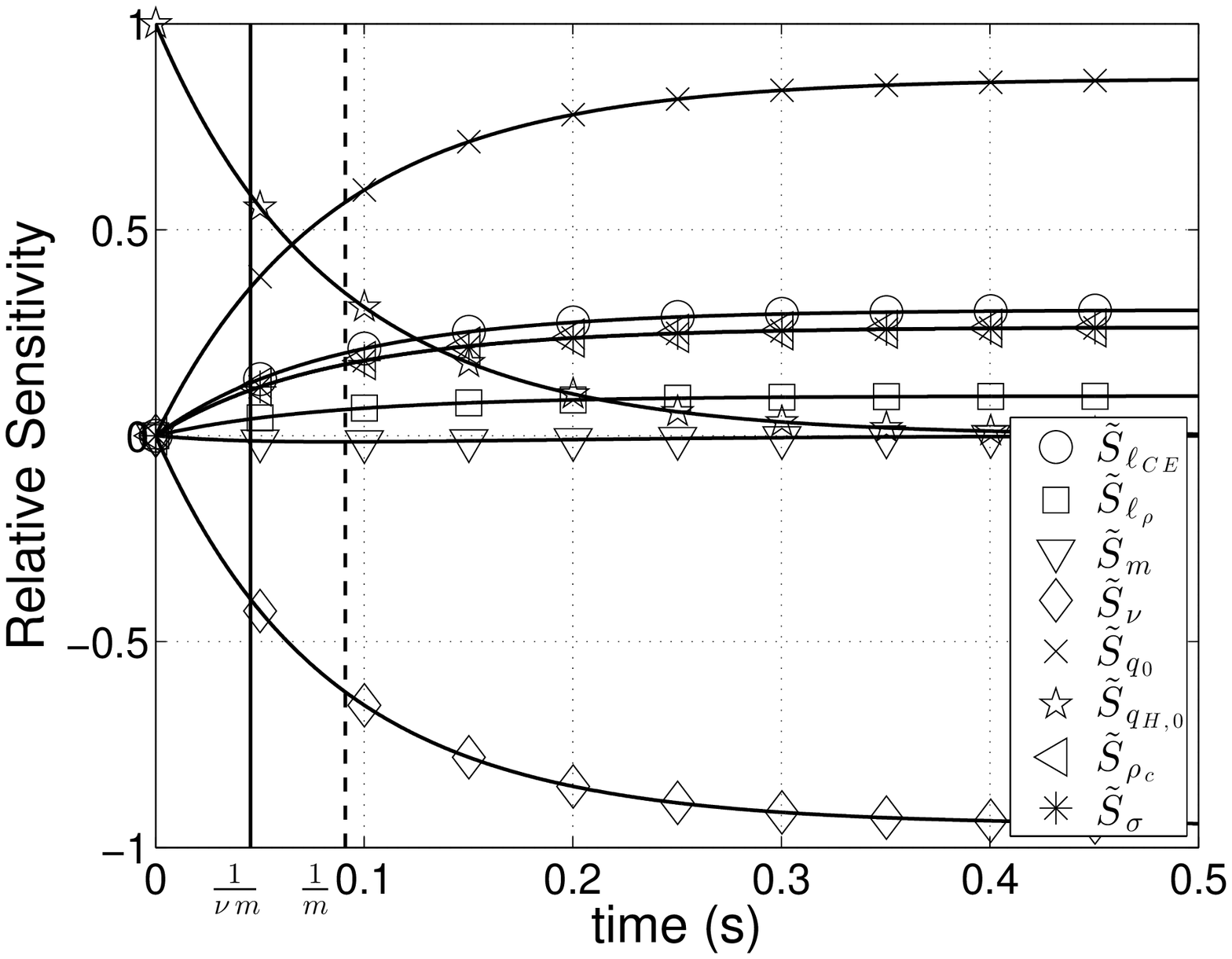} & \includegraphics[width=6cm]{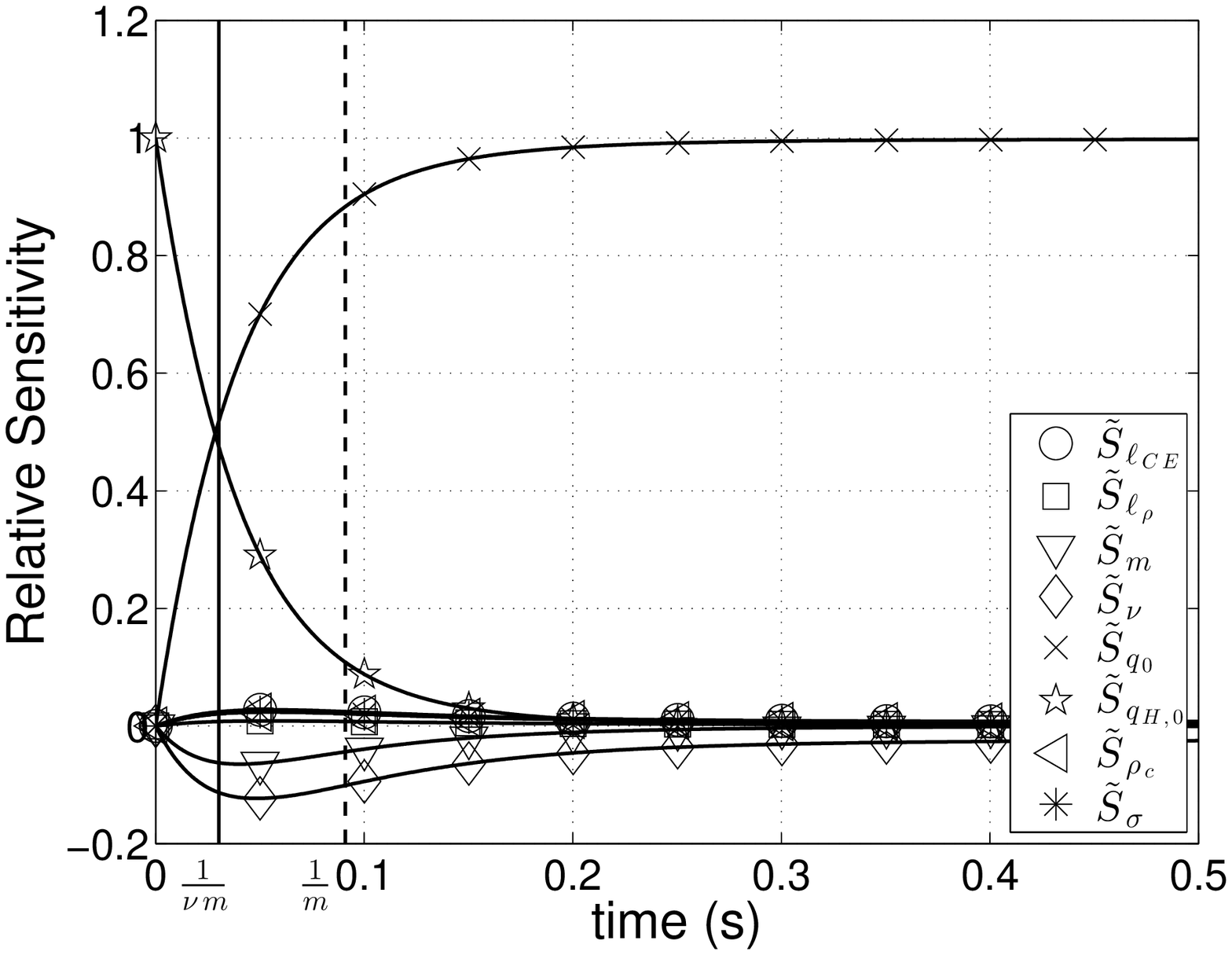} \\
\includegraphics[width=6cm]{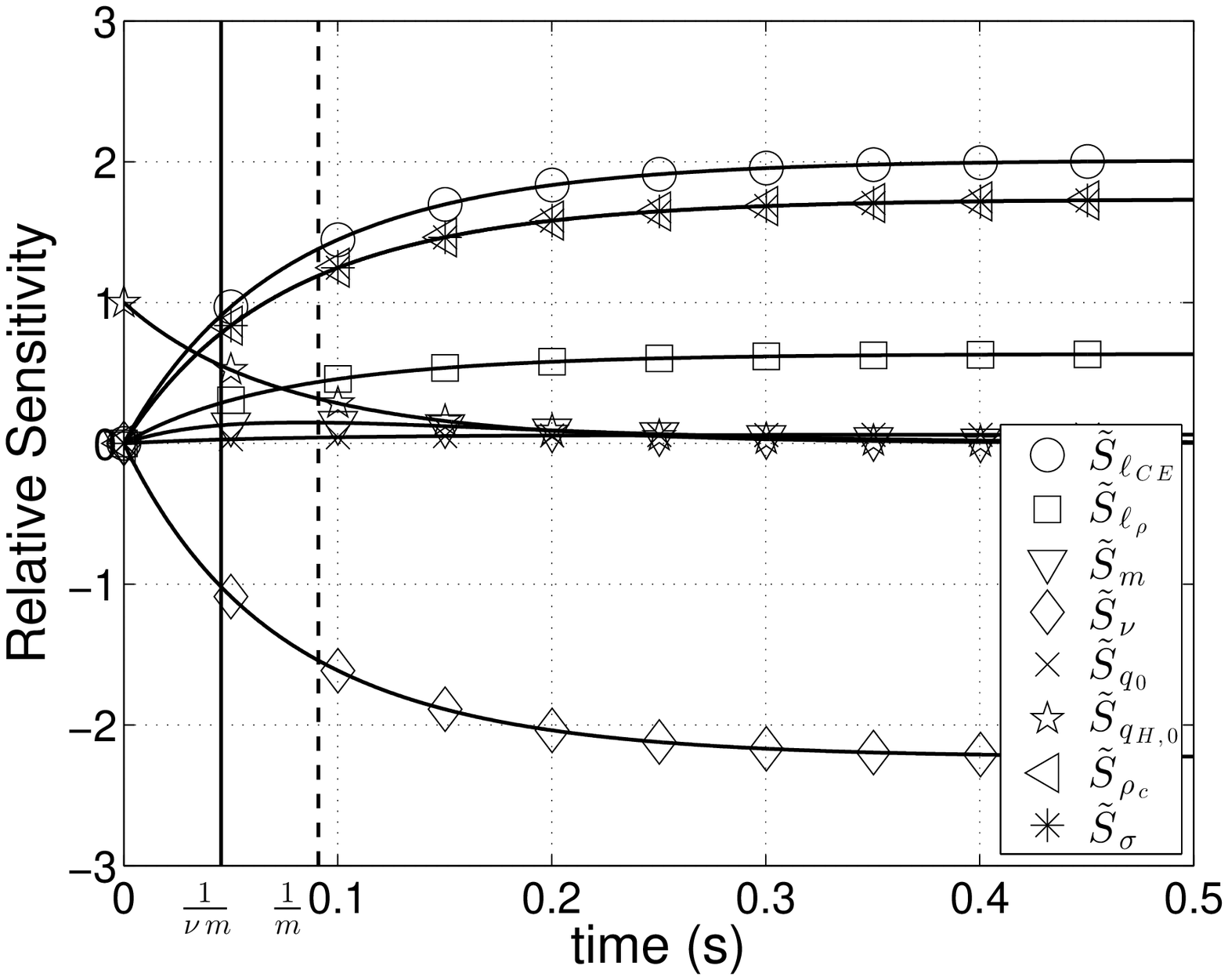} &\includegraphics[width=6cm]{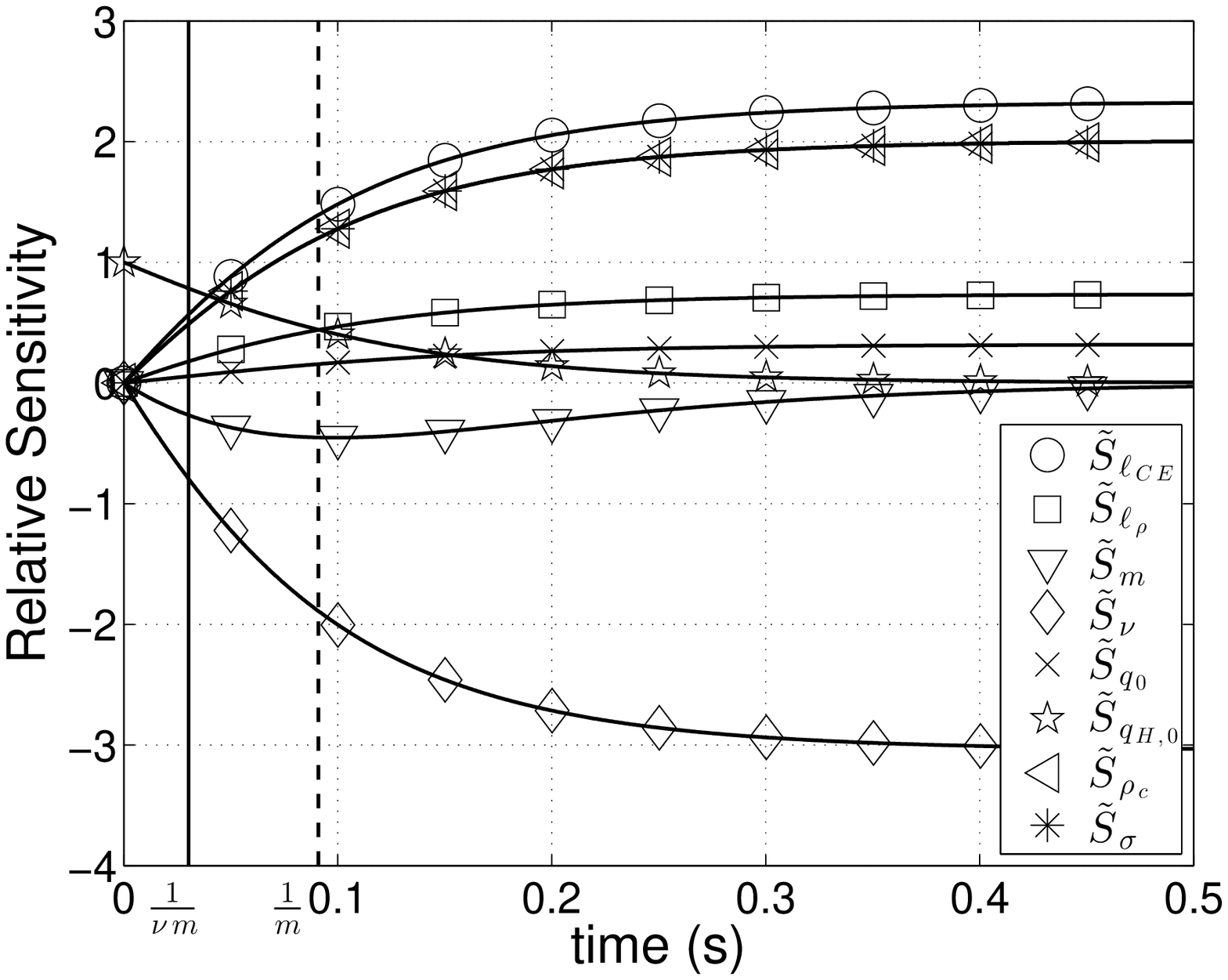} \\
\includegraphics[width=6cm]{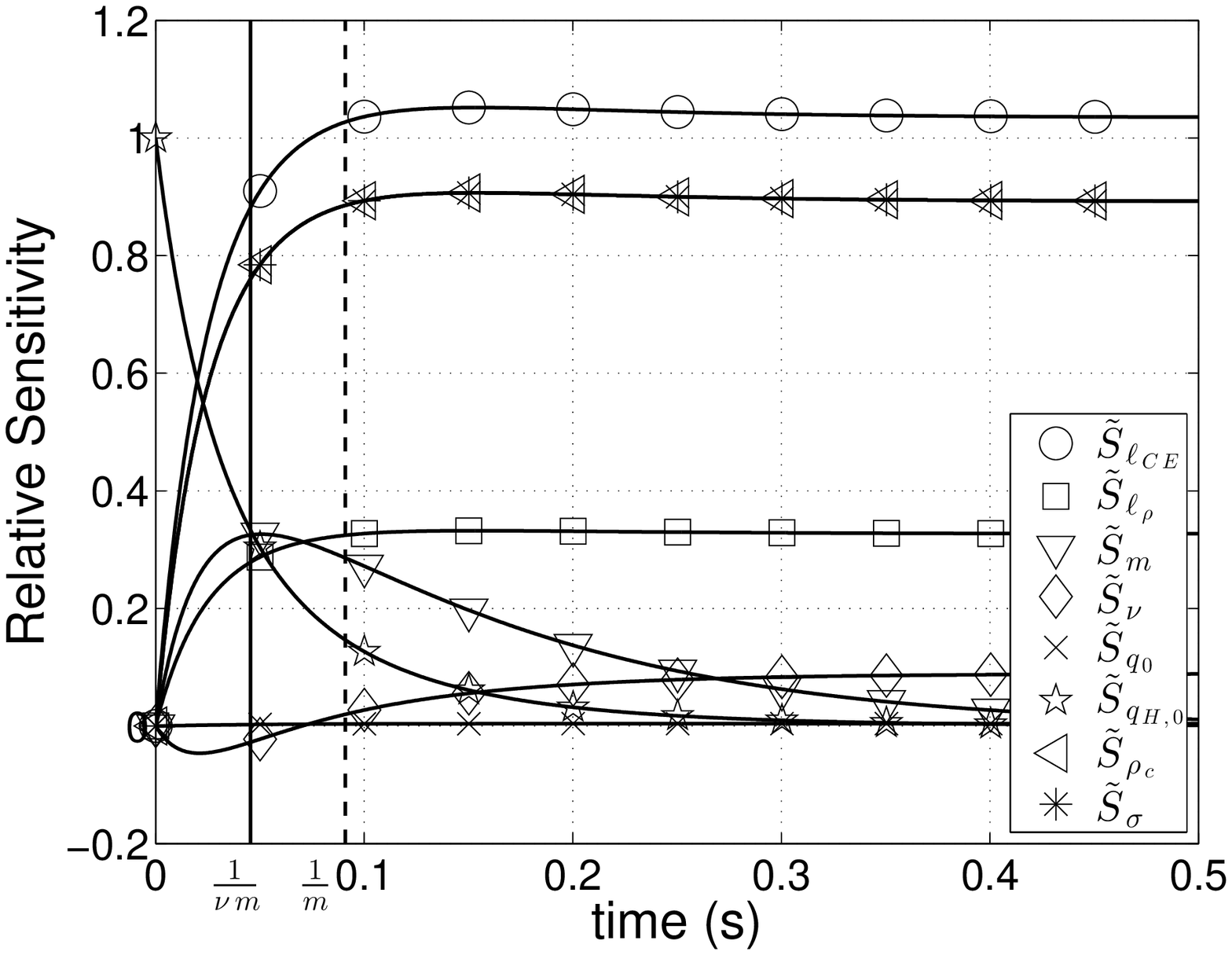} & \includegraphics[width=6cm]{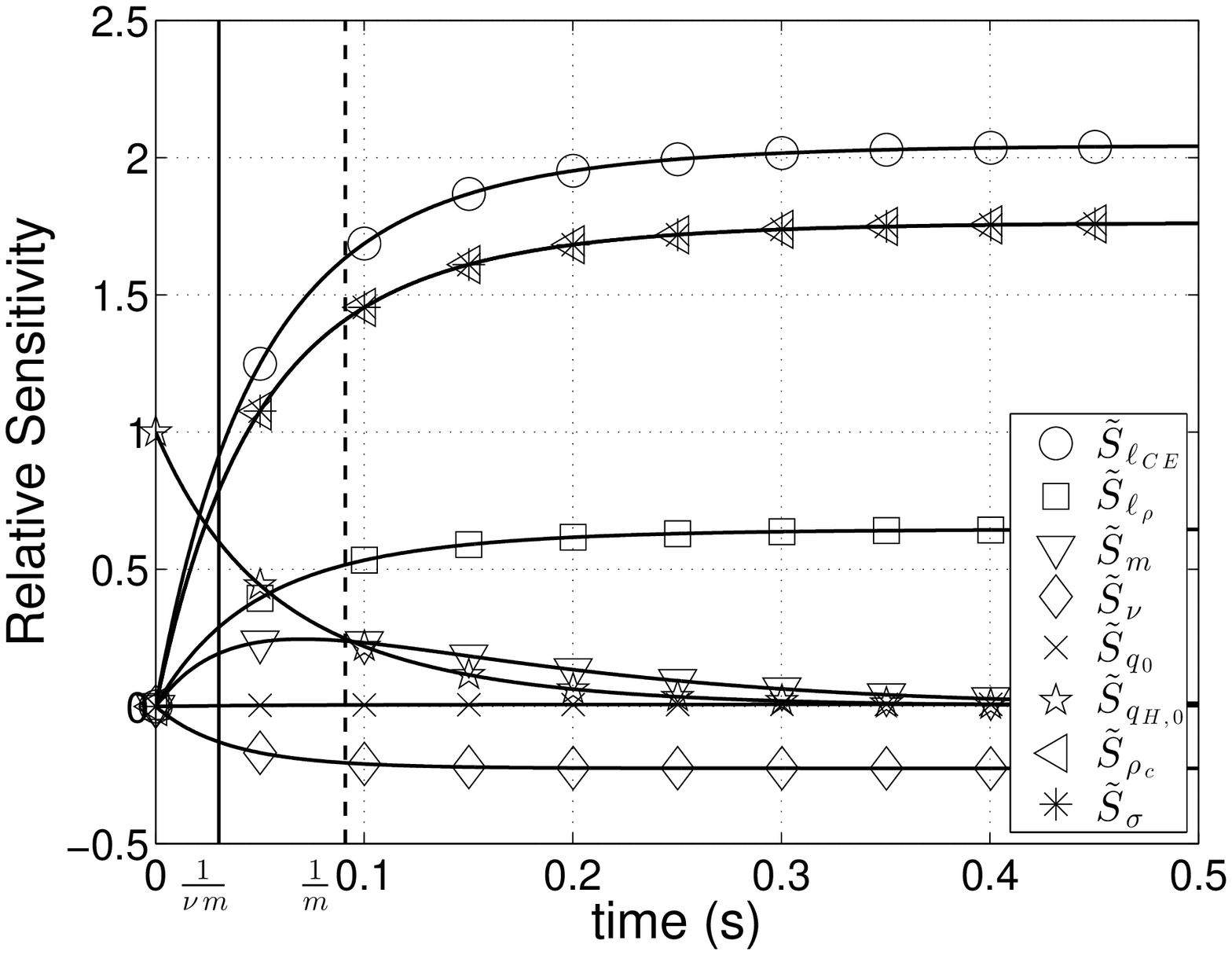} \\
\includegraphics[width=6cm]{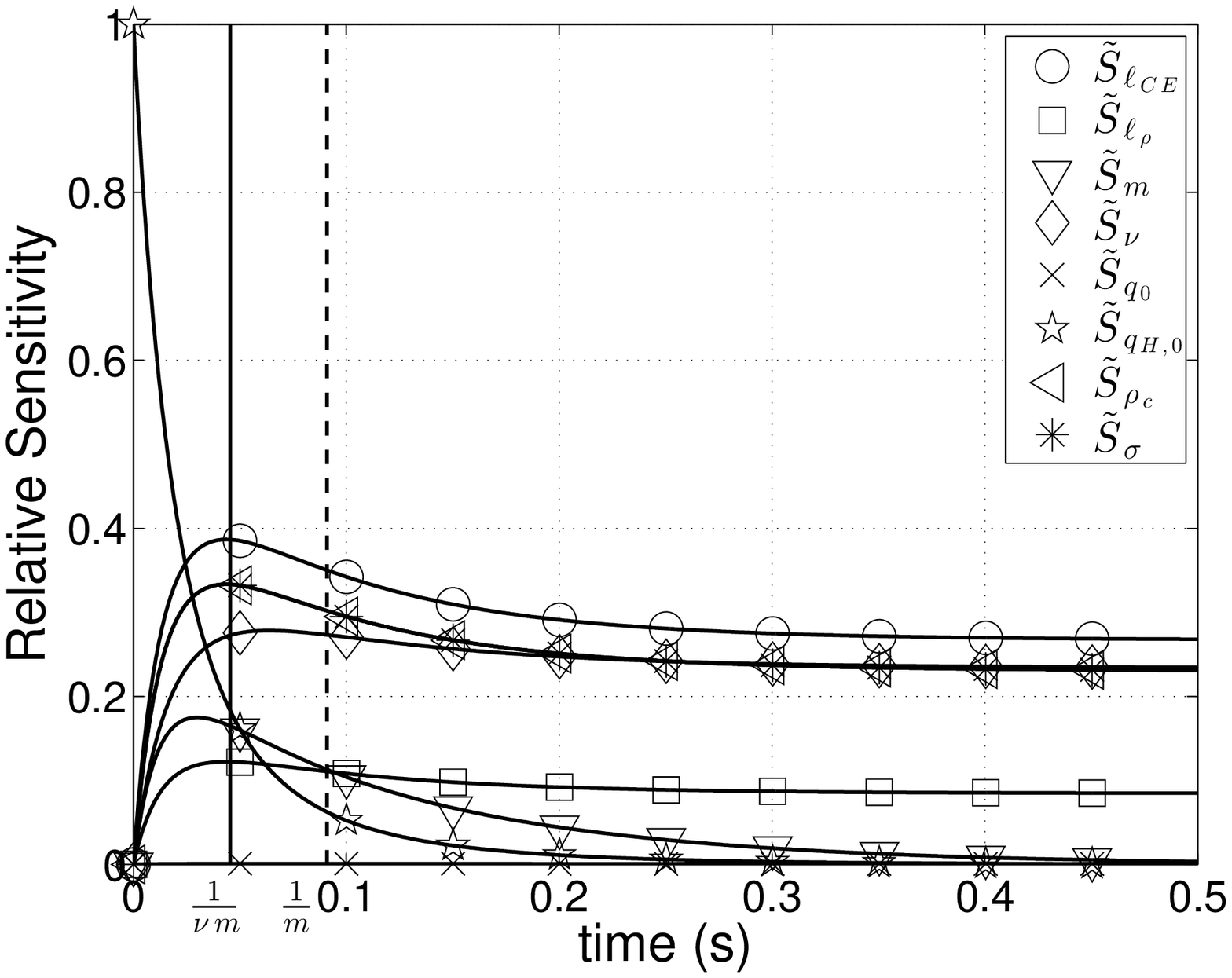} & \includegraphics[width=6cm]{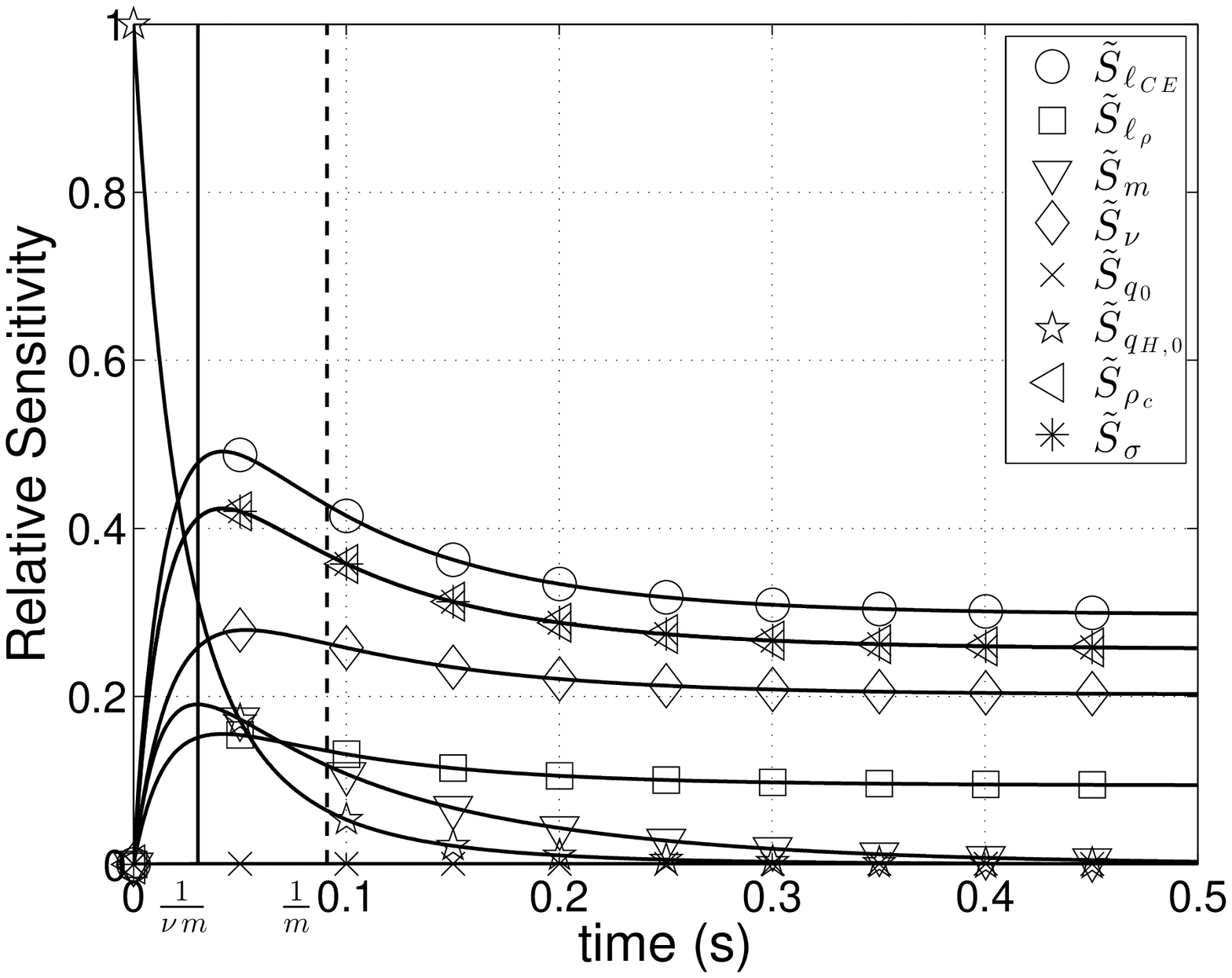} 
\end{matrix}$
\caption{Relative sensitivities $\tilde{S}_{i}$ w.r.t. all parameters
  $\lambda_i$ (set $\Lambda_H$ \eqref{lam_H}) in Hatze's activation
  dynamics \eqref{dq_H}. Parameter values varied from top (i) to
  bottom (iv) row:
\mbox{(i) $q_{H,0}=q_0=0.005,\,\sigma = 0.01$}, 
\mbox{(ii) $q_{H,0}=0.05,\,\sigma = 0.1$},
\mbox{(iii) $q_{H,0}=0.2,\,\sigma = 0.4$},
\mbox{(iv) $q_{H,0}=0.5,\,\sigma = 1$}; left column:
$\nu=2,\,\rho_c=9.10$, right column:
$\nu=3,\,\rho_c=7.24$. 
\label{fig_Hatze_result}}
\end{figure}

The latter sensitivity w.r.t. $\nu$ itself is extraordinarily high at
low activities and stimulations around $0.1$, both for $\nu=2$ and
$\nu=3$ (Fig.\,\ref{fig_Hatze_result}, second row from top), additionally
at extremely low levels for $\nu=2$
(Fig.\,\ref{fig_Hatze_result}, left, top row). At moderately submaximal
levels (Fig.\,\ref{fig_Hatze_result}, third row from top), the
solution is then influenced with an already inverted tendency
($\tilde{S}_{\nu}$ changes sign to positive) after around an $1/m$
time horizon for $\nu=2$. Here however, the solution is practically
insensitive to $\nu$ for any $\nu$. At high levels
(Fig.\,\ref{fig_Hatze_result}, bottom row) then, we find that
there is no change in the character of time
evolution of the solution, whatever specific $\nu$ value. The
degree of non-linearity $\nu$ does not matter because the time evolution
and the ranking of all other sensitivities is hardly influenced by
$\nu$. In both cases anyway, the rise in activity is speeded up by
increasing $\nu$ ($\tilde{S}_{\nu}>0$), as opposed to low activity and
stimulation levels where rises in activity are slowed down
($\tilde{S}_{\nu}<0$; see also above).
\subsubsection*{Relative sensitivities $\tilde{S}_{\sigma}$, $\tilde{S}_{\rho_c}$,
$\tilde{S}_{\ell_{CErel}}$, $\tilde{S}_{\ell_{\rho}}$}
Of all the remaining parameters, stimulation $\sigma$, scaled maximum
free Ca$^{2+}$-ion concentration $\rho_c$, relative CE length
$\ell_{CErel}$, and the pole $\ell_{\rho}$ of the length dependency in
Hatze's activation dynamics, the latter has the lowest influence on
the solution, whereas the influence characters of these four
parameters are completely identical. That
is, their sensitivities are always positive and coupled by fixed
scaling ratios due to all of them occurring within just one product
on the right side of \eqref{dq_H}. $\tilde{S}_{\sigma}$ and
$\tilde{S}_{\rho_c}$ are identical, and the sensitivity
w.r.t. $\ell_{CErel}$ is the highest, with a ratio
$\tilde{S}_{\ell_{CErel}}/\tilde{S}_{\ell_{\rho}} \approx 3$
and $\tilde{S}_{\ell_{CErel}}/\tilde{S}_{\sigma} \approx 1.2$. Except
at very low activity, where $q_0$ plays a dominating role,
and except for the generally changeful $\nu$ influence, these
are the four parameters that dominate the solution after an initial
phase in which the initial activity $q_{H,0}$ determines its
evolution. The parameter $m$ does not have a strong direct influence
on the solution. As said above, it yet defines the approximate time
horizon in which the $q_{H,0}$ influence gets lost and all other
parameters' influence is switched on from zero at $t=0$.
\clearpage
\subsection{Variance-based sensitivity and total sensitivity indices for Zajac's and Hatze's activation
  dynamics}
\label{Total_sensitivity_result}

In the table below we give the lower and upper boundaries for every parameter in $\Lambda_Z$ 
and $\Lambda_H$ used in our calculations. We do not relate any individual boundary value to a 
literature source but refer to \citet{Hatze1}, \citet{Zajac} or \citet{GS} for traceability of our choices.

\begin{table}[h]
\begin{tabular}{c ||c|c|c|c|c|c|c|c|c|c}
Parameter &  $\beta$ &  $\ell_{CErel}$ & $\ell_p$ &  $m$ & $\nu$ & $q_0$
 &  $q_{Z,0} ,q_{H,0}$ & $\rho_c$ & $\sigma$ & $\tau$ \\\hline
Lower bound & $0.1$ & $0.4$ & $2.2$ & 3 & 1.5 &0.001  & 0.01  & 4 &  0 & 0.01 \\\hline
Upper bound & 1 & $1.6$ & $3.6$  & 11 & 4 & 0.05 & 1 & 11 & 1  & 0.05    \\
\end{tabular}  
\label{bounds}   
\end{table}

\begin{figure}[!ht]
\includegraphics[width=6.5cm]{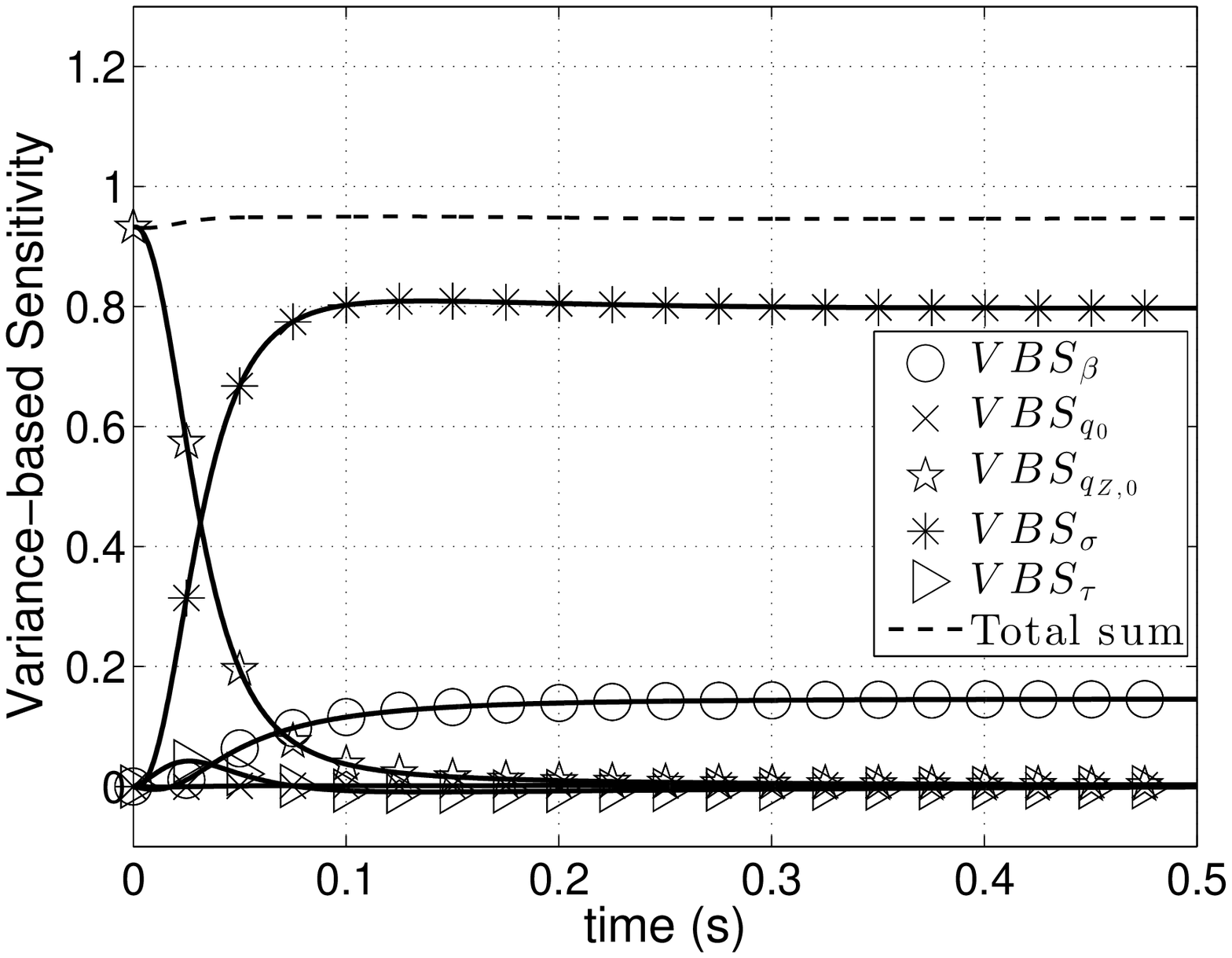}  \includegraphics[width=6.5cm]{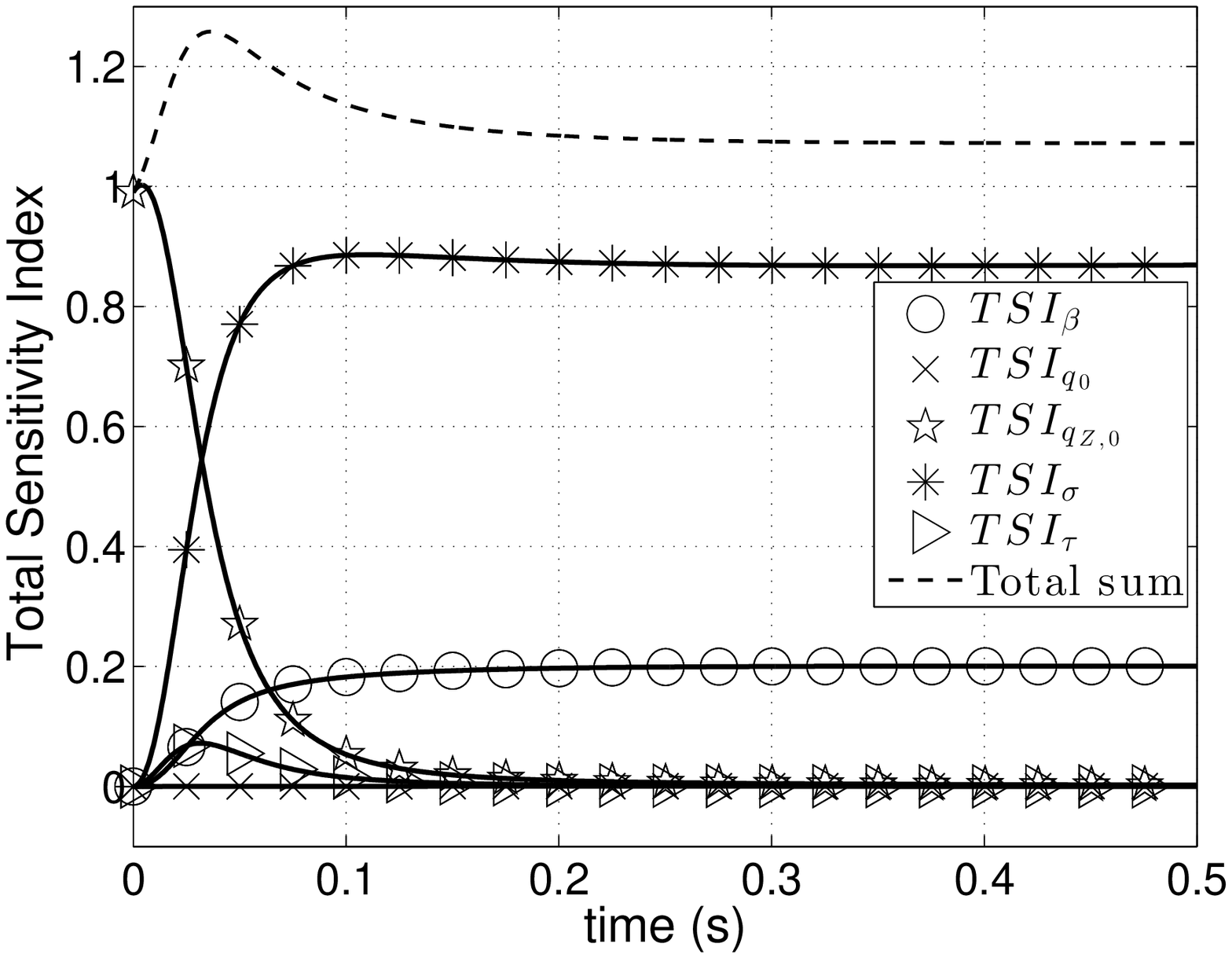}  
\caption{Variance-based sensitivity (left) and total sensitivity index (right) of every
parameter of Zajac's activation dynamics equation}
 \label{VBS_Zajac}
\end{figure}

The left hand side of Fig.~\ref{VBS_Zajac} shows the variance-based sensitivity functions 
of every parameter in $\Lambda_Z$ of Zajac's model. We compare these results to our 
previously computed relative first order sensitivity functions from Fig.~\ref{fig_Zajac_result}: 
At first sight $\tilde{S}_{q_{Z,0}}$ and $VBS_{q_{Z,0}}$ look equal but the variance based 
sensitivity function increases the duration of influence of $q_{Z,0}$ a little. For $\tau$ the 
$VBS$ also peaks at the typical time from $\tilde{S}_\tau$ but with a smaller amplitude. 
The behaviour of $VBS_\sigma$ and $VBS_\beta$ is also comparable to $\tilde{S}_\sigma$ 
and $\tilde{S}_\beta$ from the second and third row of Fig.~\ref{fig_Zajac_result}. 
Additionally we plotted the sum of all first order sensitivities. This sum indicates which 
amount of the total variance is covered by first order sensitivities. The closer the sum 
is to 1 the less second and higher order sensitivities occur.

On the right hand side of Fig.~\ref{VBS_Zajac} we see the total sensitivity index functions 
of every parameter in $\Lambda_Z$ of Zajac's model. It is noticeable that the $TSI_i$ look 
similar to the previous mentioned $VBS_i$. We interpret the graphic results as an 
importance of every parameter as suggested by \citet{Chan1997a}. Hence, the importance 
of $q_{Z,0}$ is only measurable at the beginning of the activation. At $t=0$ the importance 
is nearly $100\%$ but exponentially vanishing. The parameter $\tau$ is just of a little 
importance while the activation build-up is in progress. After saturating at a constant 
activation level the importance is only shared between $\sigma$ (major importance) 
and $\beta$ (minor importance).

\begin{figure}[!ht] 
\includegraphics[width=6.5cm]{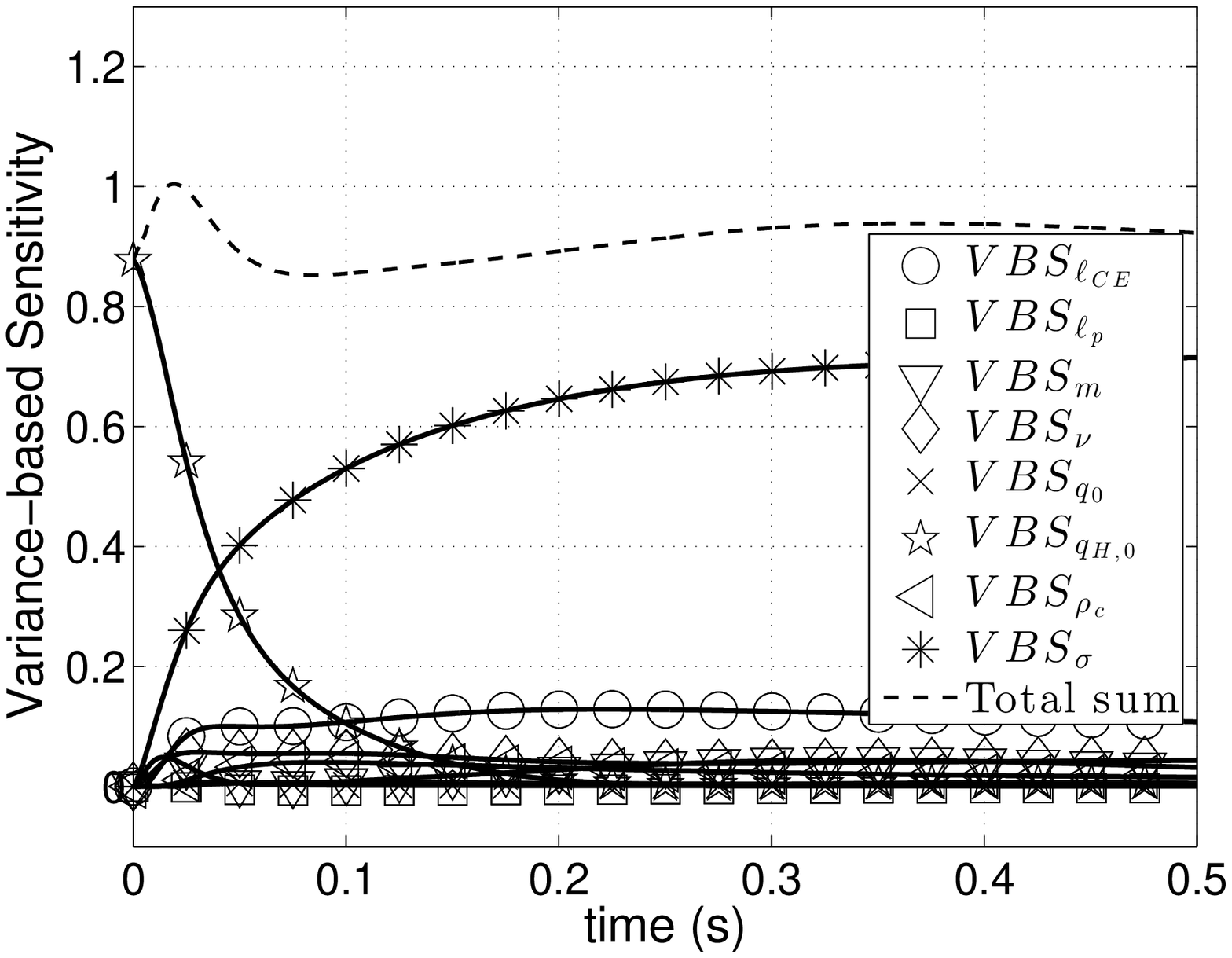} \includegraphics[width=6.5cm]{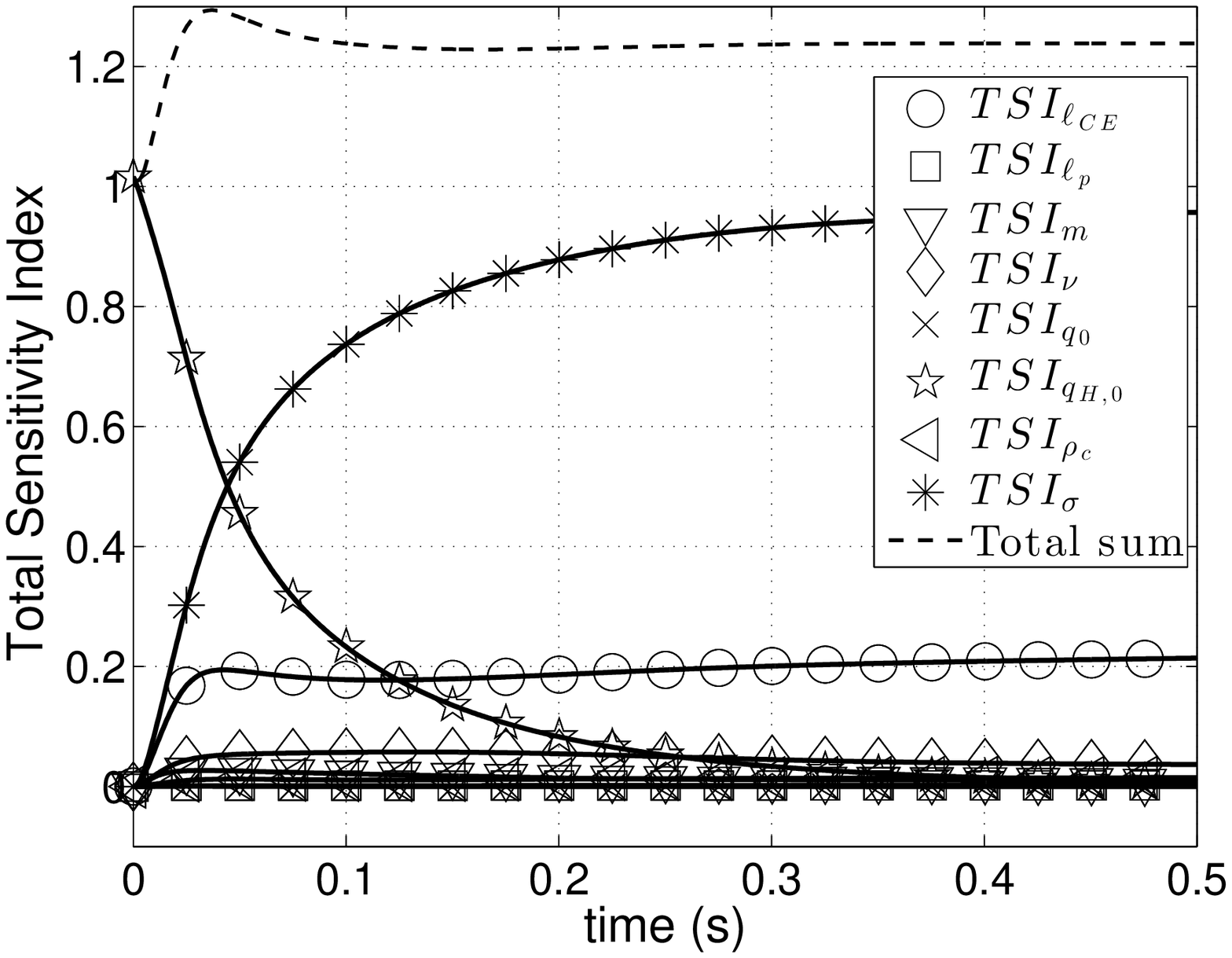} 
\caption{Variance-based sensitivity (left) and total sensitivity index (right) 
of every parameter of Hatze's activation dynamics equation}
\label{VBS_Hatze}
\end{figure}

The left hand side of Fig.~\ref{VBS_Hatze} shows the variance-based sensitivity functions 
of every parameter in $\Lambda_H$ of Hatze's model. The curve shape of $VBS_{q_{H,0}}$ 
is similar to $VBS_{q_{Z,0}}$ and $\tilde{S}_{q_{H,0}}$. For $\sigma$ the $VBS$ is again 
comparable to the second and third row in Fig.~\ref{fig_Hatze_result}. $VBS_m$ is peaking 
in a small value and similar to $VBS_\tau$. The main differences are $VBS_{\ell_{CErel}}$, 
$VBS_\nu$, $VBS_{\ell_p}$ and $VBS_{\rho_c}$ which are significant lower than the respective 
relative sensitivity functions.

 On the right hand side we see again the total sensitivity indices of Hatze's model. As above 
 the $TSI_i$ and $VBS_i$ look alike and allow an interpretation in the meaning of importance. 
 Hence, $q_{H,0}$ is as important for Hatze's model as is $q_{Z,0}$ for Zajac's. In the steady 
 state of Hatze's dynamics, which is reached after a longer time period than in Zajac's, the 
 importance is again split almost exclusively between $\sigma$ (major importance) and 
 $\ell_{CErel}$ (minor importance). All other parameter are of almost negligible importance.

\clearpage
\section{Consequences, discussion, conclusions} \label{discuss}
\subsection{A bottom line for comparing Zajac's and Hatze's activation
  dynamics: second order sensitivities} \label{firstresume}
At first sight, Zajac's activation dynamics \citep{Zajac} is more
transparent because it is descriptive in a sense that it captures the
physiological behaviour of activity rise and fall in an apparently
simple way. It thereto utilises a linear differential equation with
well-known properties, allowing for a closed-form
solution. It needs only a minimum number of parameters to describe the
Ca$^{2+}$-ion influx to the muscle as a response to electrical
stimulation: the stimulation $\sigma$ itself as a control parameter,
the time constant $\tau$ for an exponential response to a step
increase in stimulation, and a
third parameter $\beta$ (deactivation boost) biasing both the rise
time and saturation value of activity depending on stimulation and
activity levels. The smaller $\beta<1$ (deactivation in fact slowed
down compared to activation), the faster is the very
activity level $q_{Z}\,|_{_{\beta=1}} = q_{0} + \sigma\cdot (1-q_{0})$
reached, at which saturation would occur for $\beta=1$. Saturation for
$\beta<1$ occurs at a level
$q_{Z} = q_{0} + (1-q_{0})/((1-\beta)+\beta/\sigma)$ that is higher than
$q_{Z}\,|_{_{\beta=1}}$. Altogether, in Zajac's as
compared to Hatze's activation dynamics, the outcome of setting a control
parameter value $\sigma$ in terms of how fast and at which level the
activity saturates seems easier to be handled by a controller.

A worse controllability of Hatze's activation dynamics
\citep{Hatze} may be expected from its
non-linearity, a higher number of parameters, and their
interdendent influence on model dynamics. Additionally, Hatze's
formulation depends
on the CE length $\ell_{CErel}$, which makes the mutual coupling of
activation with contraction dynamics more interwoven. So, at first
sight, this seems to be a more intransparent construct for a
controller to deal with a muscle as the biological actuator. Regarding
the non-linearity exponent $\nu$, solution sensitivity further depends
non-monotonously on activity level, partly even with the strongest
influence, partly without any influence. We also found that the
solution is more sensitive to its parameters $\sigma$, $\ell_{CErel}$,
$\ell_{\rho}$ than is Zajac's activation dynamics to any of its
parameters.

This higher complexity of Hatze's dynamics becomes even more evident by
analysing the second order sensitivities (see \eqref{S2} as
well as
\eqref{S2norm} for their relative values). They express how a first
order sensitivity changes upon variation of any other model
parameter. In other words, they are a measure of model
entanglement and complexity. Here, we found that the highest values amongst
all relative second order sensitivities in Zajac's activation dynamics
are about $-0.8$ ($\tilde{R}_{\beta\,\sigma}$) and $1.6$
($\tilde{R}_{\beta\,\beta}$). In Hatze's activation dynamics, the
highest relative second order sensitivities are those with respect to
$\nu$ or $\ell_{CErel}$ (in particular for $\sigma$, $\rho_c$ and
$\nu$, $\ell_{CErel}$
themselves) with maximum values between about
$-8.0$ ($\tilde{R}_{\ell_{CErel}\,\nu}$, $\tilde{R}_{\nu\,\rho_c}$) and
$13.4$ ($\tilde{R}_{\ell_{CErel}\,\ell_{CErel}}$,
$\tilde{R}_{\ell_{CErel}\,\rho_c}$, $\tilde{R}_{\ell_{CErel}\,\sigma}$,
$\tilde{R}_{\nu\,\nu}$ at submaximal activity). That is, they 
are an order of magnitude higher than in Zajac's activation dynamics.

Yet, we have to acknowledge that Hatze's activation dynamics contains
crucial physiological features that go beyond Zajac's description.
\subsection{A plus for Hatze's approach: length dependency}
\label{lengthdep}
It has been established that the length dependency of activation
dynamics is both physiological
\citep{Kistemaker2005a} and functionally vital \citep{Kistemaker2007b}
because it largely contributes to low-frequency muscle stiffness. It
was also verified that Hatze's model approach provides a
good approximation for experimental data \citep{Kistemaker2005a}. In that
study, $\nu=3$ was used without comparing to the $\nu=2$ case. There
seem to be arguments in favour of $\nu=2$ from a mathematical point of
view. Especially, the less changeful scaling of the activation
dynamics' characteristics down to very low activity and stimulation
levels, and particularly a remaining CE length sensitivity of the
dynamics, seem to be an advantage when compared to the $\nu=3$
case. This applies subject to being in accordance with physiological
reality. It seems that experimental data with a good resolution of
activation dynamics as a response to very low muscular stimulation
levels are missing in literature so far: the lowest analysed level in
\citet{Kistemaker2005a} was $\sigma=0.08$, i.e., comparable to the
second rows from top in
Figs.\,\ref{fig_Zajac_result},\ref{fig_Hatze_result}.
\subsection{An optimal parameter set for Hatze's
  activation dynamics plus CE force-length relation}
\label{newoptimise}
Sensitivity analysis allows to rate Hatze's approach as an
entangled construct. Additionally, \citet{Kistemaker2005a}
decided to choose $\nu=3$ without giving a reason for discarding
$\nu=2$. Also, it seemed that they did not perform an algorithmic
optimisation on their muscle parameters to fit known shifts in optimal
CE length $\Delta \ell_{CE,isom,max}$ at submaximal stimulation levels,
i.e., the CE length value where the submaximal isometric force
$F_{isom} = F_{isom}(q, \ell_{CE})$ peaks. Accordingly, it seemed
worth to perform such an optimisation because $F_{isom}$ generally
depends on length $\ell_{CE}$ and activity $q$, and the latter may be
additionally biased by an $\ell_{CE}$-dependent capability for
building up cross-bridges at a given level $\gamma$ of free
Ca$^{2+}$-ions in the sarcoplasma, as formulated in Hatze's approach:
$F_{isom}(q, \ell_{CE}) = F_{max} \cdot q(\gamma,\ell_{CE}) \cdot F_{\ell}(\ell_{CE})$.
Thus, a shift in optimal CE length $\Delta \ell_{CE,isom,max}$ with
changing $\gamma$ can occur depending on the specific choices of both
the length-dependency of activation $q(\gamma,\ell_{CE})$ (see
\eqref{q_H},\eqref{rho}) and the CE's force-length relation
$F_{\ell}(\ell_{CE})$.

Consequently, we searched for optimal parameter sets of Hatze's
activation dynamics in combination with two different
force-length relations $F_{\ell}(\ell_{CE})$:
either a parabola \citep{Kistemaker2005a}
or bell-shaped curves \citep{GS,Moerl2012a}. For a given
optimal CE length $\ell_{CE,opt}=14.8\,mm$ \citep{Siebert2014a}
representing a rat gastrocnemius muscle and three fixed exponent
values $\nu=2,3,4$ in Hatze's activation dynamics (all other
parameters as given in section \ref{Model}), we thus determined
Hatze's constant $\rho_{0}$ and the width parameters of the two
different force-length relations $F_{\ell}(\ell_{CE})$ ($WIDTH$ in
\citet{Kistemaker2005a,vanSoest1993b} and $\Delta W_{asc}=\Delta
W_{des}=\Delta W$ in \citet{Moerl2012a}, respectively) by an
optimisation approach. The objective function to be minimised was the
sum of squared differences between the $\Delta \ell_{CE,isom,max}$
values as predicted by the model and as derived from experiments (see
Table 2 in \citet{Kistemaker2005a}) over five stimulation levels
$\sigma=0.55,0.28,0.22,0.17,0.08$. Note that $\gamma=\sigma$
applies in the isometric situation (see \eqref{Gamma} and compare
\eqref{q_H}).

The optimisation results are summarised in Table \ref{Optim}. The
higher the $\nu$ value, the smaller is the
optimisation error. Along with that decrease the predicted width
values $WIDTH$ or $\Delta W$, respectively. We would, however, tend to
exclude the case $\nu=4$ because the predicted width values seem
unrealistically low when compared to published values from other
sources (e.g., $WIDTH=0.56$ \citep{vanSoest1993b}, $\Delta W=0.35$
\citep{Moerl2012a}).
Furthermore, $\rho_{0}$ decreases with $\nu$ using the parabola model
for $F_{\ell}(\ell_{CE})$ whereas it saturates between $\nu=3$ and
$\nu=4$ for the bell-shaped model. The
bell-shaped model shows the most realistic $\Delta W$ in the case
$\nu=3$ ($\Delta W=0.32$). Fitting the same model to other contraction
modes of the muscle \citep{Moerl2012a}, a value of $\Delta W=0.32$ had
been found . In contrast, when using the parabola model, realistic
$WIDTH$ values between $0.5$ and $0.6$ are predicted by our
optimisation for $\nu=2$. When
comparing the optimised parameter values across all start values of
the $F_{\ell}(\ell_{CE})$ widths, across all $\nu$ values, and
across both $F_{\ell}(\ell_{CE})$ model functions, we find that the
resulting optimal parameter sets are more consistent for
bell-shaped $F_{\ell}(\ell_{CE})$ than for the parabola function.
The bell-shaped force-length relation gives generally a better
fit. For each single $\nu$ value, the corresponding optimisation error
is smaller when comparing realistic, published $WIDTH$ and $\Delta W$
values that may correspond to each other ($WIDTH=0.56$
\citep{vanSoest1993b} and $\Delta W=0.35$
\citep{Moerl2012a}). Additionally, the error values from our
optimisations are generally smaller than the corresponding value
calculated from Table 2 in \citet{Kistemaker2005a} ($0.23\,mm$).

In a nutshell, we would say that the most realistic model for the
isometric force $F_{isom}$ at submaximal activity levels is the
combination of Hatze's approach for activation dynamics with $\nu=3$
and a bell-shaped curve for the force-length relation
$F_{\ell}(\ell_{CE})$ with
$\nu_{asc}=3$. As a side effect, we predict that the parameter value
$\rho_{0}$, being a weighting factor of the first addend in the compact
formulation of Hatze's activation dynamics \eqref{dq_H}, should be
reduced by about 40\% ($\rho_{0} = 3.25 \cdot 10^{4}\,\frac{l}{mol}$)
as compared to the value originally published in \citet{Hatze2}
($\rho_{0} = 5.27 \cdot 10^{4}\,\frac{l}{mol}$).

\subsection{A generalised method for calculating parameter sensitivities}
\label{genaralmethod}
The findings in the last section were initiated by thoroughly
comparing two different biomechanical models of muscular activation
using a systematic sensitivity analysis as introduced in
\citet{Dickinson} and \citet{Lehman1982a}, respectively. Starting with
the latter formulation, \citet{Scovil} calculated specific parameter
sensitivities for muscular contractions. They applied three variants
of this method:

Method 1 applies to state variables that are explicitly known to
the modeller as, for example, an eye model \citep{Lehman1982a},
a musculo-skeletal model for running that includes a Hill-type muscle
model \citep{Scovil}, or the activation models analysed in our
study. \citet{Scovil} calculated the change in the value of a state
variable averaged over time per a finite change in a parameter value,
both normalised to each their unperturbed values. They thus calculated
just one (mean) sensitivity value for a finite time interval (e.g., a
running cycle) rather than time-continuous sensitivity functions.

Method 2: Whereas \citet{Dickinson} and \citet{Lehman1982a} had
introduced the full approach for calculating such sensitivity
functions, \citet{Scovil} distorted this approach by suggesting that the
partial derivative of the right hand side of an ODE, i.e., of the
{\em rate of change} of a state variable, w.r.t. a model parameter
would be a ``model sensitivity''. The distortion becomes explicitly
obvious from our formulation: this partial derivative is just one of
two addends that contribute to the rate of change of the sensitivity
function \eqref{S}, rather than it defines the sensitivity of the
state variable itself (i.e., the solution of the ODE) w.r.t. a model
parameter \eqref{FO}.

Method 3: \citet{Scovil} had also asked for calculating the influence
of, for example, a parameter of the activation dynamics (like the time
constant) on an arbitrary joint angle, i.e., a variable
that quantifies the overall output of a coupled dynamical system. Of
course, the time constant does not explicitly appear in the mechanical
differential equation for the acceleration of this very joint angle,
which renders applicability of method 2 impossible. The conclusion in
\citet{Scovil} was to apply method 1.  Here, the
potential of our formulation comes particularly to the fore. It enables to
calculate the time-continuous sensitivity of all components of the
coupled solution, i.e., any state variable $y_{k}(t)$. This is because
all effects of a parameter change are in principle reflected within
{\em any} single state variable, and the time evolution of a
sensitivity according to \eqref{S} takes this into account.

In this paper, we have further worked out the sensitivity function approach
by \citet{Lehman1982a}, presenting the differential equations for
sensitivity functions in more detail to those modellers who want 
to apply the method. Furthermore, we enhanced the approach by
\citet{Lehman1982a} to also calculating the sensitivities of the state 
variables w.r.t. their initial conditions \eqref{y_0}. This should be
helpful not only in biomechanics but also, for example, in 
meteorology when predicting the behaviour of storms
\citep{Langland}. Since initial conditions are often just known
approximately but start with the relative sensitivity values of $1$,
their influence should be traced to verify how their uncertainty
propagates during a simulation. In the case of muscle activation
dynamics, the sensitivities $\tilde{S}_{q_{Z,0}}$ and
$\tilde{S}_{q_{H,0}}$, respectively, rapidly decreased to zero:
initial activity has no effect on the solution early before steady
state is reached.

Furthermore, we included a second order sensitivity analysis which is
not only helpful for a deeper understanding of the parameter influence
but also part of mathematical optimisation techniques \citep{Sunar}.
The values of $\tilde{R}_{ijk}$ could be either interpreted as the
relative sensitivity of the sensitivity $\tilde{S}_{ik}$
w.r.t. another parameter $\lambda_j$ (and vice versa: $\tilde{S}_{jk}$
w.r.t $\lambda_i$) or as the curvature of the graph of the solution
$y_k(t)$ in the $N+M$-dimensional solution-parameter space. The latter
may help to connect the results to the field of mathematical
optimisation in which the second derivative (Hessian) of a function is
often included in objective functions to find optimal parameter sets.

\subsection{Deeper understanding through global methods}

As a last point for discussion we want to give some additional conclusions arising 
from the use of a global sensitivity analysis method. In section \ref{Total_sensitivity_result} 
we presented the variance-based sensitivity and the total sensitivity index according 
to \citet{Chan1997a}. In the case of Zajac's activation dynamics, we can strengthen 
our assumption that there are no significant second and higher order sensitivities 
with the exception of activation build-up. For an experimenter the only way to get 
information about the activation time constant $\tau$ is through looking at the first 
few milliseconds after a change in stimulation, but keeping in mind that the influence 
is biased by the other parameters.

For Hatze's activation dynamics, we see that there are higher order sensitivities even in 
the steady state case. When we talked about controllability of the models we presumed 
that Zajac's dynamics would be easier to handle than Hatze's. But Fig.~\ref{VBS_Hatze} 
shows that the stimulation is as well the most important control factor with even a higher 
importance than in Zajac's formulation. 

Another, at first sight unapparent, result is the importance of $\rho_c$. From a strictly 
differential analysis we concluded that this parameter should have the same sensitivity 
as $\sigma$ since they both are linear factors in Hatze's ODE which holds true for their 
relative sensitivities. The importance, in contrast, is significantly smaller, almost negligible. 
An explanation can be found if we look at the respective ranges. In the product 
$\rho_c\cdot \sigma \in [4;11]\times[0;1]$ the parameter $\rho_c$ has a wider range 
but only serves as a lever for $\sigma$ which has a much larger percentage changeability. 
The same observation can be made for the parameter $\nu$ which has a very small 
variability throughout literature. Although the differential sensitivity is quite large, $\nu$ 
has not much importance for the model output. 

Nevertheless the advantages of a global sensitivity analysis the findings must be treated 
with caution because we have a whole dynamical system summed up to a single function 
per parameter. In sections \ref{Zajac_result} and \ref{Hatze_result} we evaluated the local 
sensitivity of Zajac's and Hatze's formulation each at 8 different points in the parameter 
space. Therefore we saw the effects of each parameter on the solution in some borderline 
cases which are averaged in a global analysis. 

Summarizing the findings of this article it takes three components to a meaningful sensitivity analysis:
a deep understanding of the model, a complete mathematical investigation and an interpretation 
of the results based on the model itself.

%\section*{References}
%\bibliographystyle{plain}
\bibliographystyle{elsarticle-harv}
\bibliography{Bib}
%
%
%\vspace{2cm}
%Hereby we assure that this work has been exclusively done by the four
%authors and that there is no conflict of interest.
%
\newpage
\begin{table}[h]
\renewcommand{\baselinestretch}{1.25}
%\tablesize
\renewcommand{\arraystretch}{1.5}
\renewcommand{\tabcolsep}{2.0ex}
\caption{\label{Optim}
Parameters minimising the sum over five submaximal stimulation levels
$\gamma=\sigma=0.55,0.28,0.22,0.17,0.08$ of squared differences between
shifts in optimal CE length $\Delta \ell_{CE,isom,max}(\gamma)$ ($\Delta
l_{MA,opt}$ by Roszek et al. (1994) in third column of Table 2 in
\citet{Kistemaker2005a}) at these levels predicted by the model with
the isometric force
$F_{isom}(q, \ell_{CE}) = F_{max} \cdot q(\gamma=\sigma,\ell_{CE}) \cdot F_{\ell}(\ell_{CE})$
and by experiments; simulated data represent a rat gastrocnemius
muscle with an optimal CE length $\ell_{CE,opt}=14.8\,mm$
\citep{Siebert2014a}; start value of $\rho_{0}$ was $6.0 \cdot
10^{4}\frac{l}{mol}$;
the exponents of the bell-shaped force-length relations
$F_{\ell}(\ell_{CE})$ were fixed according to \citet{Moerl2012a}
($\nu_{asc}=3$, $\nu_{des}=1.5$), the corresponding width values in
the ascending and descending branch were assumed to be equal:
$\Delta W_{asc}=\Delta W_{des}=\Delta W$;
\citet{vanSoest1993b} and \citet{Kistemaker2005a} used a parabola for
$F_{\ell}(\ell_{CE})$; for all other model parameters see sections
\ref{newoptimise} and \ref{Model}; optimisation was done by
$fminsearch$ (Nelder-Mead algorithm) in MATLAB with error
tolerances of $10^{-8}$; $error$ is the square-root of the above
mentioned sum divided by five; corresponding error value given in
Table 2 in \citet{Kistemaker2005a} was $0.23\,mm$.}
\hfill\\
\begin{tabular}{c|ccc|ccc}
  $\nu$ & \multicolumn{3}{c|}{bell-shaped \citep{GS,Moerl2012a}} 
  & \multicolumn{3}{c}{parabola \citep{vanSoest1993b,Kistemaker2005a}}\\
  \hline
        & \multicolumn{3}{c|}{$\Delta W_{start}=0.25$} &  \multicolumn{3}{c}{$WIDTH_{start}=0.46$} \\
        & $\Delta W\,[]$ & $\rho_{0}\,[10^{4}\frac{l}{mol}]$ & $error\,[mm]$ & $WIDTH\,[]$ 
        & $\rho_{0}\,[10^{4}\frac{l}{mol}]$ & $error\,[mm]$ \\
    2   & 0.46 & 3.80 & 0.08 & 0.63 & 8.78 & 0.10\\
    3   & 0.32 & 3.25 & 0.05 & 0.41 & 5.45 & 0.07\\
    4   & 0.26 & 3.20 & 0.02 & 0.34 & 4.60 & 0.05 \\
  \hline
        & \multicolumn{3}{c|}{$\Delta W_{start}=0.35$} &  \multicolumn{3}{c}{$WIDTH_{start}=0.56$} \\
        & $\Delta W\,[]$ & $\rho_{0}\,[10^{4}\frac{l}{mol}]$ & $error\,[mm]$ & $WIDTH\,[]$ 
        & $\rho_{0}\,[10^{4}\frac{l}{mol}]$ & $error\,[mm]$ \\
    2   & 0.45 & 3.80 & 0.07 & 0.53 & 6.92 & 0.11 \\
    3   & 0.32 & 3.30 & 0.05 & 0.41 & 5.67 & 0.07 \\
    4   & 0.26 & 3.20 & 0.02 & 0.34 & 4.55 & 0.05\\
  \hline
        & \multicolumn{3}{c|}{$\Delta W_{start}=0.45$} &  \multicolumn{3}{c}{$WIDTH_{start}=0.66$} \\
        & $\Delta W\,[]$ & $\rho_{0}\,[10^{4}\frac{l}{mol}]$ & $error\,[mm]$ & $WIDTH\,[]$ 
        & $\rho_{0}\,[10^{4}\frac{l}{mol}]$ & $error\,[mm]$ \\
    2   & 0.45 & 3.78 & 0.07 & 0.55 & 7.35 & 0.11\\
    3   & 0.32 & 3.25 & 0.05 & 0.41 & 5.35 & 0.07\\
    4   & 0.26 & 3.20 & 0.02 & 0.34 & 4.56 & 0.05 \\
\end{tabular}  
\end{table}
\end{document}